\documentclass[11pt,notitlepage]{amsart}

\newcommand{\e}{\textrm{e}}
\newcommand{\E}{\mathbb{E}}
\newcommand{\R}{\mathbb{R}}
\newcommand{\Prob}{\mathbb{P}}

\newcommand{\Ordo}{\mathcal{O}}
\newcommand{\Z}{\mathbb{Z}}
\newcommand{\N}{\mathbb{N}}

\newcommand{\udx}{\mathrm{dx}}
\newcommand{\udy}{\mathrm{dy}}
\newcommand{\udz}{\mathrm{dz}}
\newcommand{\udw}{\mathrm{dw}}

\newcommand{\ai}{\mathrm{Ai}}
\newcommand{\dai}{\mathrm{Ai'}}
\newcommand{\fdet}{\mathrm{det}}

\newcommand{\lp}{\left}
\newcommand{\rp}{\right}

\newtheorem{theorem}{Theorem}[section]
\newtheorem{lemma}{Lemma}[section]

\newenvironment{prf}{\hspace{-15pt}\textbf{Proof:} } {$\;$}

\newenvironment{axiom-T}{\hspace{-15pt}\textbf{Axiom-T} } {$\;$}  
\newenvironment{axiom-S}{\hspace{-15pt}\textbf{Axiom-S} } {$\;$} 
\newenvironment{axiom-M}{\hspace{-15pt}\textbf{Axiom-M} } {$\;$}  
\newenvironment{axiom-PH}{\hspace{-15pt}\textbf{Axiom-PH} } {$\;$} 
\newenvironment{axiom-R}{\hspace{-15pt}\textbf{Axiom-R} } {$\;$}

\usepackage{epsfig}
\usepackage{graphicx}
\usepackage{amsmath,amssymb}

\numberwithin{equation}{section}

\title[Gaussian fluctuations in the Airy process]{Local Gaussian 
fluctuations in the Airy and Discrete PNG processes}
\author[J. H\"agg]{Jonas H\"agg}
\address{Jonas H\"agg \newline \indent Royal institute of Technology}
\email{jonas.hagg@math.kth.se}

\begin{document}

\begin{abstract}
We prove that the Airy process, $\mathcal{A}(t)$, locally fluctuates like
a Brownian motion. 

In the same spirit we also show that, in a certain scaling limit, the so 
called 
Discrete polynuclear growth process (PNG) behaves like a Brownian motion.  
\end{abstract}

\maketitle

\section{Introduction}

\subsection{The Airy process}
The central object of study in this paper is the local behavior of the 
Airy process, $t \rightarrow \mathcal{A}(t)$, $t \in \R$, \cite{PS}. 
The Airy process is a one dimensional process with continuous paths,
\cite{Jo1}, \cite{PS}.
The interest in this process is mainly due to the fact that it is the 
limit of a number of processes appearing in the random matrix 
literature. One example is the top curve in Dyson's Brownian motion, 
see \cite{Dy},
which, when appropriately rescaled, converges to the Airy process, 
see for instance \cite{FNH} and \cite{Jo2}. Another 
example is the boundary of the north polar region in the Aztec diamond,
see \cite{EKLP1}, \cite{EKLP2} and \cite{Jo3}, 
a discrete process also converging to the Airy process, \cite{Jo3}. 
A third example, 
the Discrete polynuclear growth model (PNG), \cite{Jo2}, \cite{KS}, 
will be described in some detail
in section \ref{discpol} where we also state a theorem about its local 
(in a certain sense) fluctuations.  

A precise definition of $\mathcal{A}(t)$ goes as follows:

The extended Airy kernel, \cite{FNH}, \cite{M}, \cite{PS}, is defined by
\begin{equation}
A_{s,t}(x,y) = \Bigg \{
\begin{array}{rl}
\int_0^{\infty} e^{-z (s-t)} \ai(x+z) \ai(y+z) \, \udz & 
\textrm{if } s \geq t \\
-\int_{-\infty}^0 e^{z (t-s)} \ai(x+z) \ai(y+z) \, \udz & 
\textrm{if } s < t ,
\end{array}
\end{equation}
where $\ai$ is the Airy function. $A_{s,s}(x,y)$ is easily seen to be the 
ordinary Airy kernel, \cite{TW}. Given $\xi_1, \ldots , \xi_m \in \R$ and 
$t_1 < \ldots < t_m$ in $\R$ we define $f$ on 
$\{ t_1, \ldots ,t_m \} \times \R$ by
\begin{equation*}
f(t_i,x) = \chi_{(\xi_i,\infty)}(x).
\end{equation*}
It is shown in \cite{Jo2} that 
\begin{equation*}
f^{1/2}(s,x) A_{s,t}(x,y) f^{1/2}(t,y)
\end{equation*}
is the integral kernel of a trace class operator on
$L^2(\{ t_1, \ldots ,t_m \} \times \R)$ where we have counting measure
on $\{ t_1, \ldots ,t_m \}$ and Lebesgue measure on $\R$.
The Airy process, $t \rightarrow \mathcal{A}(t)$, is the stationary 
stochastic process with finite dimensional distributions given by
\begin{equation*}
\Prob 
\lp[ \mathcal{A}(t_1) \leq \xi_1, \ldots ,\mathcal{A}(t_m) \leq \xi_m \rp]
= \fdet \lp( I - f^{1/2} A f^{1/2} \rp)_{L^2(\{ t_1, \ldots ,t_m \} \times \R)}.
\end{equation*}
The determinant in the right hand side is a so called Fredholm determinant.

Our main theorem states that if we condition the Airy process to be 
at some given point at time $t_1$ it will then behave, on a local scale,
like a Brownian motion. 
\begin{theorem}\label{theorem1}
Let $\epsilon > 0$ be small, $t_1 \in \R$ and 
$t_i = t_{i-1} + s_{i} \epsilon$, $2 \leq i \leq m$, where
$s_2, \ldots , s_{m} >0$. Also, let $p_1 \in \R$ and
define the sets $A_i$, $i=2, \ldots ,m$, by
\begin{equation*}
A_i = \lp\{ x \in \R | p_1 + a_i \sqrt{\epsilon} \leq x \leq 
p_1 + b_i \sqrt{\epsilon} \rp\}
\end{equation*}
where $a_i,b_i$ are given real numbers. It holds that 
\begin{multline*}
\Prob \lp[ \mathcal{A}(t_2) \in A_2 , \ldots , \mathcal{A}(t_m) \in A_m \rp
| \mathcal{A}(t_1) = p_1] \\ =
\int_{a_2}^{b_2} dx_2 \cdots \int_{a_m}^{b_m} dx_m 
\frac{1}{\sqrt{4 \pi s_2}}  e^{-\frac{x_2^2}{4 s_2}} 
\prod_{i=3}^m \frac{1}{\sqrt{4 \pi s_i}} 
e^{-\frac{(x_i - x_{i-1})^2}{4 s_i}} + E
\end{multline*}
where 
\begin{equation*}
|E| \leq \sqrt{\epsilon} \log{\epsilon^{-1}} \prod_{i=2}^m(b_i-a_i) 
C_{p_1,s_2, \ldots ,s_m}.
\end{equation*}
\end{theorem}
Figure 1 describes the setup in the theorem.
\begin{figure}[hb]
\begin{center}
\includegraphics{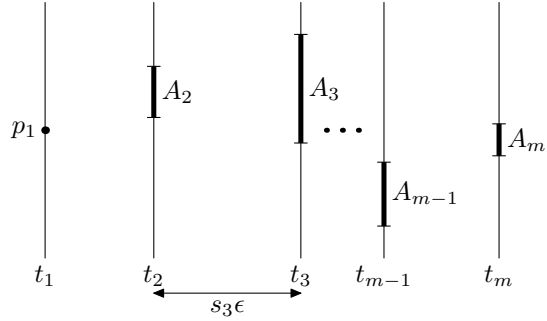}
\caption{Conditioned that $\mathcal{A}(t_1)=p_1$ Theorem \ref{theorem1}
gives the approximate probability for the process to move through 
the sets $A_i$. Note that $t_{i+1}-t_i \sim \epsilon$ and 
$|A_i| \sim \sqrt{\epsilon}$.}
\end{center}
\end{figure}

\emph{Remark 1.} A couple of previous results about the Airy process are
the following:

In \cite{PS} it is shown that 
\begin{equation*}
\mathrm{Var}(\mathcal{A}(t) - \mathcal{A}(0)) = 2t + \Ordo(t^2)
\end{equation*}
as $t \rightarrow 0$.

In \cite{AvM}, see also \cite{W}, the long distance 
covariance asymptotics for the Airy process is calculated to be
\begin{equation*}
\E \lp[ \mathcal{A}(t) \mathcal{A}(0) \rp] - 
\E \lp[ \mathcal{A}(t)\rp] \E \lp[ \mathcal{A}(0) \rp] =
t^{-2} + \Ordo(t^{-4})
\end{equation*}
as $t \rightarrow \infty$.
This proves that $\mathcal{A}(t)$ is not a Markov process since this would
imply exponential decay.
 
\emph{Remark 2.} Given Theorem \ref{theorem1} it is natural to ask the 
corresponding question about processes converging to the Airy process.
Theorem \ref{PNG} in section \ref{discpol} below provides such a result  
for the Discrete polynuclear growth process.

\subsection{The extended Airy point process}
We now present another construction, \cite{Jo2}, of the Airy process that
will help us analyzing its local behavior.

Let $m \in \Z_+$ be arbitrary and $t_1 < t_2 \ldots < t_m$ be points in $\R$ 
which we 
shall think of as times. Define
\begin{equation*}
E = \R_{t_1} \cup \R_{t_2} \cup \cdots \cup \R_{t_m}.
\end{equation*}
We shall refer to $\R_{t_j}$ as time line $t_j$. We define $X$ to be the space
of all locally finite countable configurations of points (or particles) in $E$.
Locally finite means that, if $x=(x_1,x_2, \ldots) \in X$ then, for any bounded
set $C \subset E$, it holds that $\# (C \cap x) < \infty$. Here $\# B$ 
represents the
number of points in the set $B$. One can construct a $\sigma$-algebra on $X$
from the cylinder sets: Let $B \subset E$ be any bounded Borel set and 
$n \geq 0$.
Define
\begin{equation*}
C_n^B = \lp\{ x \in X : \# B = n\rp\}
\end{equation*}
to be a cylinder set and $\Sigma$ to be the minimal 
$\sigma$-algebra that contains all cylinder sets.
One can now define probability measures on the space $(X,\Sigma)$. 
The so called
extended Airy point process is an example of such a measure and it will be 
described below.

For the sake of convenience, we will often denote the extended Airy kernel by
$A(x,y)$ instead of $A_{t_i,t_j}(x,y)$ when it is clear that $x \in \R_{t_i}$
and $y \in \R_{t_j}$. Let $z_1, \ldots , z_k$ be points in $E$. 
The $k$-point correlation function is defined by
\begin{equation} \label{corr_def}
R(z_1, \ldots , z_k) = \fdet \lp[ A(z_i,z_j) \rp]_{i,j=1}^k .
\end{equation}
It is possible to show that these correlation functions determine a 
probability measure on $(X,\Sigma)$, the extended Airy point process, by 
demanding that the following identity holds, \cite{So}:
\begin{equation} \label{corr_eq}
\E \lp[ \prod_{i=1}^n \frac{\# B_i !}{(\# B_i - k_i)!} \rp] =
\int_{B_1^{k_1} \times \cdots \times B_n^{k_n}} R(z_1, \ldots ,z_k)\, \udz .
\end{equation}
Here $B_1, \ldots , B_n$ are
disjoint Borel subsets of $E$ and $k_i \in \Z_+$, $1 \leq i \leq n$, are 
such that $k_1 + \ldots + k_n = k$. 

It is possible to show that, at each time line $\R_{t_i}$, there is almost 
surely a largest particle, $\lambda(t_i)$, and 
\begin{equation} \label{airy=extairy}
(\lambda(t_1), \ldots ,\lambda(t_m)) = (\mathcal{A}(t_1), \ldots ,
\mathcal{A}(t_m))
\end{equation}
in distribution, \cite{Jo2}.
It is through this representation that we are able to show that the Airy
process behaves locally as a Brownian motion.

\subsection{Discrete polynuclear growth} \label{discpol}
The second object of interest in this paper is the 
the so called Discrete polynuclear growth model (PNG), \cite{Jo2}, \cite{KS}. 
It is defined by
\begin{equation}
h(x,t+1) = \max{(h(x-1,t),h(x,t),h(x+1,t))} + \omega(x,t+1)
\end{equation}
where $x \in \Z$, $t \in \N$, $h(x,0)=0$ $\forall x \in \Z$ and 
$\omega(x,t+1) = 0$ if $|x| > t$ or if $t-x$ is even, otherwise
$\omega(x,t+1)$ are independent geometric random variables with
\begin{equation}\label{geometric}
\Prob [ \omega(x,t+1) = m ] = (1 - q) q^m \qquad  0 < q < 1 .
\end{equation}
It is convenient to extend the process to all $x \in \R$ by setting
$h(x,t) = h(\lfloor x \rfloor,t)$.
A description of this process using words and pictures goes as follows:

At time $t=1$ a block of width one and height $\omega(0,1)$ appears over 
the interval $[0,1)$. This block then grows sideways one unit in both 
directions and at time $t=2$ two blocks of width one and heights 
$\omega(-1,2)$, $\omega(1,2)$ are placed on top of it over the intervals
$[-1,0)$ and $[1,2)$ respectively. These blocks now grow one unit in each 
direction disregarding overlaps. At time $t=3$ three new blocks are placed 
over $[-2,-1)$, $[0,1)$ and $[2,3)$. This procedure goes on producing at each
time the curve $h(x,t)$ that can be thought of as a growing interface. 
Figure 2 shows a realization for $t=1,2,3$. 

\begin{figure}[ht] \label{fig2}
\begin{center}
\includegraphics{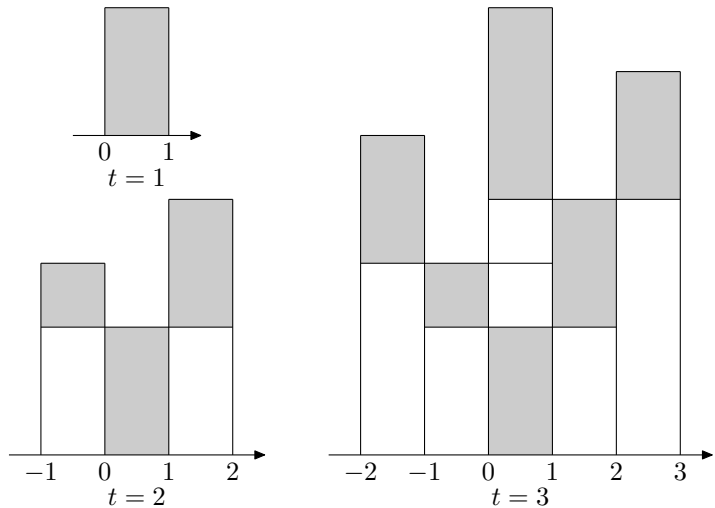}
\caption{A sample of the discrete PNG process for $t=1,2,3$. The shaded
blocks represent the growth due to the random variables $\omega (x,t)$.}
\end{center}
\end{figure}
The process $h$ is closely connected to a growth model, $G(M,N)$, studied 
in \cite{Jo1}. Let $w(i,j)$, $(i,j) \in \Z_+^2$, be independent 
random variables with distribution given by (\ref{geometric}). Define
\begin{equation*}
G(M,N) = \max_{\pi} \sum_{(i,j) \in \pi} w(i,j)
\end{equation*}
where the maximum is taken over all up/right paths from $(1,1)$ to $(M,N)$.
One can think of $G(M,N)$ as a point to point last-passage time and 
\begin{equation*}
G_{pl}(N) = \max_{|K|<N} G(N+K,N-K)
\end{equation*}
as a point to line last-passage time. In \cite{Jo2} it is shown that
\begin{equation*}
G(i,j) = h(i-j,i+j-1).
\end{equation*}
The definition of $G_{pl}$ therefore inspires the study of
$K \rightarrow h(2K,2N-1)$, that is, the height curve at even sites at time
$2N-1$.

In \cite{Jo2} the rescaled process, $t \rightarrow H_N(t)$, $t \in \R$, is, 
for appropriate $t$, defined by
\begin{equation*}
d N^{1/3} H_N(t) = 
h \lp( 2 \frac{1 + \sqrt{q}}{1 - \sqrt{q}} d^{-1} N^{2/3} t,2N-1 \rp)
- \frac{2 \sqrt{q}}{1- \sqrt{q}} N
\end{equation*}
and for the rest of $\R$ by the use of linear interpolation. The constant $d$
is given by
\begin{equation*}
d = \frac{(\sqrt{q})^{1/3}(1 + \sqrt{q})^{1/3}}{1- \sqrt{q}}.
\end{equation*}
The main result about $H_N(t)$ in \cite{Jo2}  is the following theorem:
\begin{theorem}[Johansson]\label{H_thm}
Let $\mathcal{A}(t)$ be the Airy process defined by its finite dimensional 
distributions and $T$ be an arbitrary positive number. There is a 
continuous version of $\mathcal{A}(t)$ and
\begin{equation*}
H_N(t) \rightarrow \mathcal{A}(t) - t^2
\end{equation*}
as $N \rightarrow \infty$ in the weak$^*$-topology of probability measures on 
$C(-T,T)$.
\end{theorem}
In particular this theorem shows that the fluctuations of $h$ are of order
$N^{1/3}$ and that non-trivial correlations in the transversal
direction show up when looking at times $t_i$ where
$t_{i+1} - t_i \sim N^{2/3}$.

Motivated by Theorems \ref{theorem1} and \ref{H_thm} one 
could guess that 
$h$, on a time scale of order $N^{\gamma}$, $0< \gamma < 2/3$, behaves like
a Brownian motion. The theorem below shows that this is indeed the case.

Given some $m \in \Z_+$ set 
\begin{align*}
K_1 & =  \frac{1 + \sqrt{q}}{1 - \sqrt{q}}d^{-1}N^{2/3}\tau_1 \\
K_{i+1} & = K_i + \frac{1 + \sqrt{q}}{1 - \sqrt{q}}d^{-1}s_{i+1}N^{\gamma}
\qquad i=1, \ldots ,m-1
\end{align*}
where $0 < \gamma < \frac{2}{3}$ and $\tau_1$, $s_i > 0$ are real numbers
such that $K_i \in \Z$. Define 
\begin{equation*}
J_1 = \frac{2 \sqrt{q}}{1 - \sqrt{q}}N + \psi d N^{1/3} \in \Z_+
\end{equation*}
where $\psi$ is any real number such that $J_1 \in \Z$. 
\begin{theorem}\label{PNG}
Define the sets $A_i$, $i=2, \ldots ,m$, by
\begin{equation*}
A_i = \lp\{ j \in Z_+ | 
j = J_1 + x_i d N^{\gamma / 2}, \, a_i \leq x_i \leq b_i \rp\}
\end{equation*}
where $a_i,b_i$ are given real numbers. There exists $c>0$ such that
\begin{align*}
\Prob \big[ h(2K_2, 2N-1) \in A_2 , \ldots , h(2K_m, 2N-1) \in & A_m  \\
& | \, h(2K_1, 2N-1) = J_1 \big] \\ =
\int_{a_2}^{b_2} dx_2 \cdots \int_{a_m}^{b_m} dx_m \frac{1}{\sqrt{4 \pi s_2}} 
e^{-\frac{x^2}{4 s_2}}
\prod_{i=3}^m \frac{1}{\sqrt{4 \pi s_i}} & e^{-\frac{(x_i - x_{i-1})^2}{4 s_i}}
+ E 
\end{align*}
where 
\begin{equation*}
|E| \leq N^{-c} \prod_{i=2}^m(b_i-a_i) C_{\psi_1,s_2, \ldots ,s_m} .
\end{equation*}
\end{theorem}

\section{Proof of theorem \ref{theorem1}}

The connection (\ref{airy=extairy}) shows that we can prove the theorem by 
studying the largest particle in the extended 
Airy point process at times $t_1, \ldots , t_m$.

The appearance of $C$ in formulae below should be interpreted as follows: 
There exists a positive constant which may depend on $p_i$, $s_i$, 
$i=2, \ldots , m$, 
validating the inequality to the left when inserted instead of $C$. Other
error terms will typically also depend on $p_i$, $s_i$.

Set $J_1 = [p_1 - \delta_1, p_1] \subset \R_{t_1}$ and 
$J_i = [p_i - \sqrt{\epsilon} \delta_i,p_i] \subset \R_{t_i}$,
$2 \leq i \leq m$, where $\delta_i > 0$ and 
$p_i = p_{i-1} + y_i \sqrt{\epsilon}$, $y_i \in \R$.
We also set $I_i = ( p_i,\infty)$, $i=1, \ldots , m$.

We will show that
\begin{multline}\label{target}
\lim_{\scriptstyle{\delta_1 , \ldots , \delta_m \rightarrow 0^+}}
\frac{1}{\delta_2 \cdots \delta_m}
\frac{\Prob \lp[ \# J_1 \geq 1, \ldots, \# J_m \geq 1,
\# I_1 = \ldots = I_m = 0 \rp]}
{\Prob \lp[\# J_1 \geq 1, \# I_1 = 0 \rp]} \\
= \frac{1}{\sqrt{(4 \pi)^{m-1} s_2 \cdots s_{m}}} 
e^{-\frac{y_2^2}{4s_2} - \ldots -\frac{y_m^2}{4s_m}} 
+ \mathcal{O}(\sqrt{\epsilon} \log{\epsilon}),
\end{multline}
implying Theorem \ref{theorem1}.

The first step is to show that the probabilities in the numerator and 
denominator above can be approximated by appropriate expected values. 

For $k,n \in \Z_+$ we shall use the common notation
\begin{equation*}
n^{[k]} = n (n-1) \cdots (n-k+1).
\end{equation*}
Let $J$ be an interval on some
time line and $\chi_A$ be the indicator function for the event $A$. 
Since
\begin{eqnarray*}
\#J - \chi_{\{ \#J \geq 1\}} & = &  \Bigg\{ 
\begin{array}{ll} 
k-1 & \quad \, ; \, \#J = k \geq 2 \\
0 &   \quad \, ; \, \#J = 0,1
\end{array} \\
\#J^{[2]} = \#J(\#J - 1) & = & \Bigg\{
\begin{array}{ll} 
k(k-1) & ; \, \#J = k \geq 2 \\
0 &  ; \, \#J = 0,1
\end{array}
\end{eqnarray*}
it holds that
\begin{equation} \label{ineq1}
0 \leq \#J - \chi_{\{ \#J \geq 1\}} \leq \#J^{[2]}.
\end{equation}
This together with the following facts will be useful:
\begin{multline}\label{prob_approx}
\Prob \lp[ \# J_1 \geq 1, \ldots, \# J_m \geq 1, \# I_1 = 0 \rp] \\
- \Prob \lp[ \# J_1 \geq 1, \ldots, \# J_m \geq 1,
\# I_1 = \ldots = \# I_m = 0 \rp]  =
\end{multline}
\begin{multline*}
= \Prob \lp[ \# J_1 \geq 1, \ldots, \# J_m \geq 1, \#I_1 = 0, 
(\# I_2 = \ldots = \# I_m = 0)^c \rp] \\ 
= \Prob \lp[ \# J_1 \geq 1, \ldots, \# J_m \geq 1, \#I_1 = 0, 
\cup_{i=2}^m \{ \# I_i \neq 0 \} \rp] \\ 
\leq \sum_{i=2}^m 
\Prob \lp[ \# J_1 \geq 1, \ldots, \# J_m \geq 1, \#I_i \neq \#I_1 \rp] 
\end{multline*}
We now express the probabilities in terms of expected values. If we set
\begin{equation}
T(J_i) = \# J_i - \chi_{\{ \# J_i \geq 1 \}} .
\end{equation}
then
\begin{align*}
\Prob \big[ \# J_1 \geq 1, & \ldots, \# J_m \geq 1, I_1 = 0 \big] \\
& = \E \lp[ ( \# J_1 - T(J_1)) \cdots ( \# J_m - T(J_m)) 
\cdot \chi_{\{ \#I_1 = 0 \}} \rp] \\
& = \E \lp[ (\# J_1 \cdots \# J_m + U(J_1, \ldots , J_m)) 
\cdot \chi_{\{ \#I_1 = 0 \}} \rp]
\end{align*}
where $U$ is defined by the last equality.
In view of (\ref{ineq1}) and (\ref{corr_eq}) we get, for example,
\begin{multline*}
\E \lp[ T(J_1) \cdot \# J_2 \cdots \# J_m \rp] \leq 
\E \lp[ \# J_1^{[2]} \cdot \# J_2 \cdots \# J_m \rp] \\
= \int_{J_1^2 \times J_2 \cdots \times J_m} R (x_1,x_2,\ldots,x_{m+1}) \udx  = 
\mathcal{O} (\delta_1^2 \cdot \delta_2 \cdots \delta_m).
\end{multline*}

Since $U(J_1, \ldots , J_m)$ is a sum of terms like this one (at least one 
$T(J_i)$) we see that 
\begin{equation*}
\lim_{\scriptstyle{\delta_1 , \ldots , \delta_m \rightarrow 0^+}}
\frac{1}{\delta_1 \cdots \delta_m} 
\E \lp[ U(J_1, \ldots , J_m) \cdot \chi_{\{ \#I_1 = 0 \}} \rp]  =  0.
\end{equation*}
Repetition of this argument shows together with (\ref{prob_approx}) that
\begin{align*}
\lim_{\scriptstyle{\delta_1, \ldots ,\delta_m \rightarrow 0^+}} &
\frac{1}{\delta_1 \cdots \delta_m} \Prob \lp[ \# J_1 \geq 1, 
\ldots, \# J_m \geq 1,
\# I_1 = \ldots = I_m = 0 \rp]  \\
& = \lim_{\scriptstyle{\delta_1 , \ldots , \delta_m \rightarrow 0^+}}
\frac{1}{\delta_1 \cdots \delta_m} 
\E \lp[ \# J_1 \cdots \# J_m \cdot \chi_{\{I_1 = 0\}} \rp] \\
& \, \, \, \, \, \, 
+ \mathcal{O} \lp(\sum_{i=2}^m \lim_{\scriptstyle{\delta_1, \ldots , 
\delta_m \rightarrow 0^+}} \frac{1}{\delta_1 \cdots \delta_m}
\E \lp[ \# J_1 \cdots \# J_m \cdot \chi_{\{\#I_i \neq \#I_1\}} \rp] \rp) 
\end{align*}
and also that
\begin{equation*}
\lim_{\scriptstyle{\delta_1 \rightarrow 0^+}}
\frac{1}{\delta_1} \Prob \lp[ \# J_1 \geq 1,
\# I_1 = 0 \rp] = 
 \lim_{\scriptstyle{\delta_1 \rightarrow 0^+}}
\frac{1}{\delta_1} \E \lp[ \# J_1 \cdot \chi_{\{I_1 = 0\}} \rp] .
\end{equation*}
Later it will be shown that 
\begin{equation}\label{uggly_eq}
\lim_{\scriptstyle{\delta_1 , \ldots ,\delta_m \rightarrow 0^+}}
\frac{1}{\delta_1 \cdots \delta_m}
\E \lp[ \# J_1 \cdots \# J_m \cdot \chi_{\{\#I_i \neq \#I_1\}} \rp] = 
\mathcal{O}(\sqrt{\epsilon} \log{\epsilon})
\end{equation}
but let us first be constructive. 

We want to show that
\begin{multline} \label{constr_eq}
\lim_{\scriptstyle{\delta_1 , \ldots , \delta_m \rightarrow 0^+}}
\frac{1}{\delta_1 \cdots \delta_m} 
\E \lp[ \# J_1 \cdots \# J_m \cdot \chi_{\{I_1 = 0\}} \rp] \\
= \lim_{\scriptstyle{\delta_1 \rightarrow 0^+}}
\frac{1}{\delta_1} \E \lp[ \# J_1 \cdot \chi_{\{I_1 = 0\}} \rp] \\  
\times \frac{1}{\sqrt{(4 \pi)^{m-1} s_2 \cdots s_{m}}} 
\, e^{-\frac{y_2^2}{4s_2} - \ldots -\frac{y_m^2}{4s_m}}
+ \mathcal{O}(\sqrt{\epsilon}) .
\end{multline}
To start with we need to find a representation of the left hand side 
of (\ref{constr_eq}) that is suitable for analysis. 
\begin{align*}
\E \big[ \# J_1 \cdots \# J_m \cdot \chi_{\{I_1 = 0\}} \big] 
& = \E \lp[ \# J_1 \cdots \# J_m \cdot 
\lim_{\lambda \rightarrow \infty}e^{- \lambda \#I_1} \rp] \\
& = \E \lp[ \# J_1 \cdots \# J_m \cdot \lim_{\lambda \rightarrow \infty}
\sum_{k=0}^{\infty} \frac{(e^{-\lambda} - 1)^k}{k!} I_1^{[k]} \rp] \\
&= \E \lp[ \# J_1 \cdots \# J_m \cdot
\sum_{k=0}^{\infty} \frac{(- 1)^k}{k!} I_1^{[k]} \rp] 
\\ & = \sum_{k=0}^{\infty} \frac{(-1)^k}{k!} 
\E \lp[ \# J_1 \cdots \# J_m \cdot \# I_1^{[k]} \rp]
\end{align*}
In the second equality we have used the formula
\begin{equation} \label{comb_formula}
e^{\lambda n} = \sum_{k=0}^{\infty} \frac{({e}^{\lambda} -1)^k}{k!} n^{[k]}.
\end{equation}
In the fourth equality we take the sum out of the expectation. 
By Fubini's theorem we are allowed to do this since
\begin{multline*}
\E \lp[ \# J_1 \cdots \# J_m \cdot 
\sum_{k=0}^{\infty} \frac{\# I_1^{[k]}}{k!}  \rp] 
\leq \E \lp[ \# J_1 \cdots \# J_m \cdot 
\sum_{k=0}^{\infty} \frac{\# I_1^{k}}{k!} \rp] \\
= \E \lp[ \# J_1 \cdots \# J_m \cdot e^{\# I_1} \rp] \leq 
\E \lp[ \# J_1^2 \cdots \# J_m^2 \rp]^{1/2} \,
\E \lp[ e^{2 \# I_1} \rp]^{1/2} < \infty.
\end{multline*}
In fact $\E \lp[ z^{\# I_1} \rp]$ is an entire function in $z$, \cite{So}.

Another technical issue we need to deal with is to prove that
\begin{align*}
\lim_{\scriptstyle{\delta_1, \ldots ,\delta_m \rightarrow 0^+}} &
\frac{1}{\delta_1 \cdots \delta_m} 
\sum_{k=0}^{\infty} \frac{(-1)^k}{k!} 
\E \lp[ \# J_1 \cdots \# J_m \cdot \# I_1^{[k]} \rp] \\
& = \sum_{k=0}^{\infty} \frac{(-1)^k}{k!} 
\lim_{\scriptstyle{\delta_1, \ldots ,\delta_m \rightarrow 0^+}}
\frac{1}{\delta_1 \cdots \delta_m}
\E \lp[ \# J_1 \cdots \# J_m \cdot \# I_1^{[k]} \rp] \\
& = \sum_{k=0}^{\infty} \frac{(-1)^k}{k!} 
\int_{I_1^k} (\sqrt{\epsilon})^{m-1} R (p_1,\ldots ,p_m, x_1,\ldots,x_k) \, \udx.
\end{align*}
Please recall definition (\ref{corr_def}) and note that the second equality 
is immediate from (\ref{corr_eq}).
Define $G_k(z_1, \ldots , z_m)$, $z_i \in \R_{t_i}$, by 
\begin{equation}
G_k(z_1, \ldots , z_m) = \frac{(-1)^k}{k!} 
\int_{I_1^k} R (z_1,\ldots ,z_m, x_1,\ldots,x_k) \, \udx.
\end{equation}
The identity sought for is
\begin{multline} \label{sought}
\lim_{\scriptstyle{\delta_1, \ldots ,\delta_m \rightarrow 0^+}} 
\frac{1}{\delta_1 \cdots \delta_m}
\sum_{k=0}^{\infty} 
\int_{J_1 \times \cdots \times J_m} G_k(z_1, \ldots, z_m) \, \udz \\ =
\sum_{k=0}^{\infty} (\sqrt{\epsilon})^{m-1} G_k(p_1, \ldots, p_m) .
\end{multline}
This will hold if for some neighbourhood 
$\Omega$ of $(p_1, \ldots, p_m)$ there exist constants $C_k > 0$ such that
\begin{equation*}
|G_k(z_1, \ldots , z_m)| \leq C_k
\end{equation*}
if $(z_1, \ldots , z_m) \in \Omega$ and 
\begin{equation*}
\sum_{k=0}^{\infty} C_k < \infty.
\end{equation*}
That this is indeed the case follows from calculations similar to the 
ones appearing in the proof of Lemma \ref{princ_lemma} which is given
at the end of this section. 

The following lemma can be found in \cite{Ok}:
\begin{lemma} \label{phi_theorem}
Let $\alpha > 0$, then
\begin{equation*}
\int_{-\infty}^{\infty} e^{\alpha z} \ai(x+z) \ai(y+z) \, \udz =
\frac{1}{\sqrt{4 \pi \alpha}} e^{-\frac{(x-y)^2}{4 \alpha} -
\frac{\alpha}{2}(x+y) + \frac{\alpha^3}{12}}.
\end{equation*}
\end{lemma}
In this section we call this function $\phi_{\alpha}(x,y)$ or simply
$\phi(x,y)$ when it is clear what $\alpha$ is.
From Lemma \ref{phi_theorem} and the definition of the Airy kernel it 
follows that, for $s<t$
\begin{multline*}
A_{s,t}(x,y) = \int_{0}^{\infty} e^{z(t-s)} \ai(x+z) \ai(y+z) \, \udz
- \phi_{t-s}(x,y) \\ 
=: \widetilde{A}_{s,t}(x,y) - \phi_{t-s}(x,y).
\end{multline*}
For $s \geq t$ it is convenient to set 
$\widetilde{A}_{s,t}(x,y) = A_{s,t}(x,y)$.
\begin{lemma} \label{princ_lemma}
Suppose that $1 \leq v \leq m$, $v \in \Z$. Then, for some $C$ depending on
$p_1, \ldots , p_m$,
\begin{align} \label{D_u(k)} \nonumber
& (\sqrt{\epsilon})^{m-1} \,
\int_{I_v^k} R(p_1, \ldots , p_m, x_1, \ldots , x_k) \, \udx \\
& \qquad \qquad = (\sqrt{\epsilon})^{m-1} \phi(p_1,p_2) \phi(p_2,p_3) \cdots 
\phi(p_{m-1},p_m) \\ \nonumber
& \qquad \qquad \qquad \qquad 
\times \int_{I_1^k} R(p_1, x_1, \ldots , x_k) \, \udx
+ \sqrt{\epsilon} \, \mathcal{O} \lp( (Ck)^{\frac{k+m}{2}} \rp).
\end{align}
Furthermore, if $v \geq 2$ then
\begin{align} \label{D_{1,v}} \nonumber
& (\sqrt{\epsilon})^{m-1}  \,
\int_{I_1} \udx \int_{I_v} \udy \, 
R(p_1, \ldots , p_m, x, y) \\
& \qquad  = (\sqrt{\epsilon})^{m-1} \, 
\phi(p_1,p_2) \phi(p_2,p_3) \cdots \phi(p_{m-1},p_m) \\ \nonumber
& \qquad \qquad \times \lp( \int_{I_1^2} R(p_1,x_1,x_2) \, \udx + 
\int_{I_1} R(p_1,x) \, \udx \rp)
+ \mathcal{O}(\sqrt{\epsilon} \log{\epsilon}).
\end{align}
\end{lemma}
\noindent From (\ref{D_u(k)}) we now get (\ref{constr_eq}). 

We turn now to (\ref{uggly_eq}).
Clearly
\begin{multline*}
\E \lp[ \# J_1 \cdots \# J_m \cdot \chi_{\{\#I_i \neq \#I_1\}} \rp] 
\leq \E \lp[ \# J_1 \cdots \# J_m \cdot (\#I_i - \#I_1)^2 \rp] \\
= \E \lp[ \# J_1 \cdots \# J_m \cdot (\#I_i^{[2]} + \#I_1^{[2]}
+ \#I_i + \#I_1 -2 \#I_1 \#I_i ) \rp].
\end{multline*}
We now obtain (\ref{uggly_eq}) since 
\begin{align}\label{corr_eq_2} \nonumber
(\sqrt{\epsilon})^{m-1} & 
\bigg( \int_{I_i^2} R(p_1, \ldots , p_m,x,y) \, \udx \udy +
\int_{I_1^2} R(p_1, \ldots , p_m,x,y) \, \udx \udy  \\
& \qquad + \int_{I_i} R(p_1, \ldots , p_m,x) \, \udx +
\int_{I_1} R(p_1, \ldots , p_m,x) \, \udx  \\ \nonumber
& \qquad \qquad \qquad 
- 2 \int_{I_1 \times I_i} R(p_1, \ldots , p_m,x,y) \, \udx \udy \bigg) =
\mathcal{O} (\sqrt{\epsilon} \log{\epsilon})
\end{align}
by Lemma \ref{princ_lemma}.

To get (\ref{target}) we need one more result, namely that
\begin{equation} \label{pos}
\lim_{\delta_1 \rightarrow 0^+} \frac{1}{\delta_1} 
\E [ \# J_1 \chi_{\{ \# I_1 = 0 \}}] > 0. 
\end{equation}
Let $F_2(s)$ be the Tracy-Widom distribution function corresponding to
the largest eigenvalue in the GUE, \cite{TW}. Then 
\begin{align}\label{twprimeq} \nonumber
\lim_{\delta_1 \rightarrow 0^+} \frac{1}{\delta_1}  &
\E [ \# J_1 \chi_{\{ \# I_1 = 0 \}}] \\ & =
\lim_{\delta_1 \rightarrow 0^+} \frac{1}{\delta_1}
\sum_{k=0}^{\infty} \frac{(-1)^k}{k!} 
\int_{J_1} dx_0 \int_{I_1^k} d^k x \, \fdet (A(x_i,x_j))_{0 \leq i,j \leq k} \\
& = \sum_{k=0}^{\infty} \frac{(-1)^k}{k!} \nonumber
\int_{I_1^k} \fdet (A(x_i,x_j))_{0 \leq i,j \leq k} d^k x = F'_2(p_1)
\end{align}
where in the last row $x_0 = p_1$. The last equality can be obtained by
differentiating the corresponding equality for the distribution function 
$F_2(t)$, \cite{TW}, we omit the details here. The first equality
has been shown above and the second is a special case of
(\ref{sought}). Since $F_2'(s) > 0$ for all $s \in \R$, see \cite{TW}, 
we obtain (\ref{pos}).

What is still left is to prove Lemma \ref{princ_lemma}.

\noindent \textbf{Proof of Lemma \ref{princ_lemma}:}
We start with (\ref{D_u(k)}). For $0 \leq r \leq m-1$ and
$k \geq 1$ define $D_r(k)$ by
\begin{multline*}
D_r(k) = (\sqrt{\epsilon})^r \phi(p_1,p_2) \phi(p_2,p_3)\cdots \phi(p_r,p_{r+1}) 
\int_{I_v^k} \udx  \\ \times \lp|
\begin{array}{ccccc}
A(p_{r+1},p_1) & \sqrt{\epsilon} A(p_{r+1},p_{r+2}) & \ldots & 
\sqrt{\epsilon} A(p_{r+1},p_{m}) & A(p_{r+1},x_j) \\
\vdots & \vdots & \, & \vdots & \vdots \\
A(p_{m},p_1) & \sqrt{\epsilon} A(p_{m},p_{r+2}) & \ldots & 
\sqrt{\epsilon} A(p_{m},p_{m}) & A(p_{m},x_j) \\
A(x_i,p_1) & \sqrt{\epsilon} A(x_i,p_{r+2}) & \ldots & 
\sqrt{\epsilon} A(x_i,p_{m}) & A(x_i,x_j) 
\end{array} \rp| .
\end{multline*}
In the determinant $1 \leq i,j \leq k$ and for $r=0$ we set the
empty product in front of the integral to 1.
Please note that $D_0(k)$ is equal to the left hand side in (\ref{D_u(k)}).
We let $\widetilde{D}_r(k)$ be almost the same as $D_r(k)$. 
The only difference is
that we put in $\widetilde{A}(p_{r+1},p_{r+2})$ in position (1,2) 
in the matrix instead
of $A(p_{r+1},p_{r+2})$. By using induction we shall now prove that 
\begin{equation} \label{ind_eq}
D_0(k) = D_r(k) + \sqrt{\epsilon} \, \mathcal{O} \lp( (Ck)^{\frac{k+m}{2}} \rp)
\end{equation}
for $0 \leq r \leq m-1$. Clearly (\ref{ind_eq}) holds if $r=0$.
Suppose now that (\ref{ind_eq}) holds for some $r$ such that 
$0 \leq r \leq m-2$.
By expanding the determinant in $D_r(k)$ along the first row we see that
\begin{equation}
D_r(k) = D_{r+1}(k) + \widetilde{D}_r(k).
\end{equation}
What has to be proved is hence that
\begin{equation*}
\widetilde{D}_r(k) = \sqrt{\epsilon} \, 
\mathcal{O} \lp( (Ck)^{\frac{k+m}{2}} \rp).
\end{equation*}
To do this, Hadamard's inequality will come in handy but before we
recall this inequality we present a lemma which will be frequently
used from now on.
The proof is readily obtained from Lemma 
\ref{phi_theorem} and the standard estimates, see \cite{Ol},
\begin{align*}
|\ai(x)| & \leq C_M e^{-2 |x|^{3/2}/3} \\
|\dai(x)| & \leq C_M \sqrt{|x|} e^{-2 |x|^{3/2}/3}
\end{align*}
that hold for $x \geq - M$.
\begin{lemma}\label{K_approx}
Suppose that $s<t$ and $M > 0$. For $x,y \geq -M$ and any $\lambda > 0$ 
it holds that
\begin{align*}
|A_{t,s}(x,y)| & \leq C_{M,\lambda} e^{-\lambda(x+y)} \\
A_{t,s}(x,y) & = A_{t,t}(x,y) + \mathcal{O}(t-s) \,
e^{-\lambda(x+y)} \\
A_{s,t}(x,y) & = A_{t,t}(x,y) - (1 + \mathcal{O}(t-s))
\frac{1}{\sqrt{4 \pi (t-s)}} \, e^{-\frac{(x-y)^2}{4 (t-s)}} \\ 
& \hspace{7cm} + \mathcal{O}(t-s) \, e^{-\lambda(x+y)}. 
\end{align*}
The errors depend only on $M$ and $\lambda$. Moreover, 
\begin{equation*}
|A_{s,s}(x+\alpha,y) - A_{s,s}(x,y)| \leq \alpha \, C_{M,\lambda} e^{-\lambda (x+y)}
\end{equation*}
for all $\alpha > 0$. 
\end{lemma}
Let $B = (b_{i,j})_{1\leq i,j \leq n}$, $b_{i,j} \in \R$ be a matrix. 
Hadamard's inequality states that
\begin{equation}
|\fdet B| \leq \lp( \prod_{i=1}^n \sum_{j=1}^n b_{ji}^2 \rp)^{1/2}.
\end{equation}
Below we find upper bounds for the equivalent to $\sum_{j=1}^n b_{ji}^2$ 
in the matrix appearing in $\widetilde{D}_r(k)$. 

Column 1:
\begin{equation*}
\sum_{j=r+1}^m A^2(p_{r+1},p_1) + \sum_{j=1}^k A^2(x_j,p_1) \leq C (k+m)
\end{equation*}

Column 2:
\begin{multline*}
\epsilon \, \lp( \widetilde{A}^2(p_{r+1},p_{r+2}) +
\sum_{j=r+2}^m A^2(p_{j},p_{r+2}) + \sum_{j=1}^k A^2(x_j,p_{r+2}) \rp) \\
\leq \epsilon \, \bigg\{
\begin{array}{ll}
C (k+m) & \mathrm{if} \, \, v \geq r+2 \\
Cm + C \sum_{j=1}^k (\widetilde{A}(x_j,p_{r+2}) - \phi(x_j,p_{r+2}))^2 & 
\mathrm{if} \, \, v < r+2
\end{array} 
\end{multline*}

Columns $3, \ldots , m-r$ ($r+3 \leq i \leq m$):
\begin{equation*}
\epsilon \, \lp( 
\sum_{j=r+1}^m A^2(p_{j},p_{i}) + \sum_{j=1}^k A^2(x_j,p_i) \rp) \leq 
C (k+m)
\end{equation*}

Last $k$ columns ($1 \leq i \leq k$):
\begin{multline*}
\sum_{j=r+1}^m A^2(p_{j},x_i) + \sum_{j=1}^k A^2(x_j,x_i) \\ \leq \Bigg\{
\begin{array}{ll}
\sum_{j=r+1}^{v-1} \lp( \widetilde{A}(p_{j},x_i) - \phi(p_{j},x_i) \rp)^2 +
C k e^{-2x_i} & \textrm{if } v \geq r+2 \\
C (k+m) e^{-2x_i} & \textrm{if } v < r+2
\end{array} 
\end{multline*}
Next we multiply everything together, take the square root and then integrate.
Assume that $v<r+2$.
\begin{multline*}
\int_{I_v^k} \Bigg[ C(k+m)  \, \epsilon  
\lp( C + C \sum_{j=1}^k (\widetilde{A}(x_j,p_{r+2}) - \phi(x_j,p_{r+2}))^2 \rp)\\ 
\times (C(m+k))^{m-r-2}
\, (C(k+m))^k e^{-2(x_1 + \ldots + x_k)} \Bigg]^{1/2} \, \udx \\
\leq \sqrt{\epsilon} \, (Ck)^{\frac{k+m}{2}} \int_{I_v^k} 
e^{-(x_1 + \ldots + x_k)} \lp( 
1 + \sum_{j=1}^k (1 + \phi(x_j,p_{r+2})) \rp) \, \udx \\
\leq \sqrt{\epsilon} \, (Ck)^{\frac{k+m}{2}}
\end{multline*}
The case $v \geq r+2$ can be treated similarly.

To obtain (\ref{D_u(k)}) it remains to show that 
\begin{multline*}
\int_{I_v^k} \fdet \lp[ 
\begin{array}{cc}
A(p_m,p_1) & A(p_m,x_j) \\
A(x_i,p_1) & A(x_i,x_j)
\end{array} \rp]_{1 \leq i,j \leq k} \, \udx \\
= \int_{I_1^k} \fdet \lp[ 
\begin{array}{cc}
A(p_1,p_1) & A(p_1,x_j) \\
A(x_i,p_1) & A(x_i,x_j)
\end{array} \rp]_{1 \leq i,j \leq k} \, \udx
+ \sqrt{\epsilon} \, \mathcal{O} \lp( (Ck)^{\frac{k+m}{2}} \rp).
\end{multline*}
This is quite easily achieved using Hadamard's inequality and
Lemma \ref{K_approx}.
We do not present the details here but instead go on to 
prove (\ref{D_{1,v}}).
 
The first part of the proof will be similar to the proof of (\ref{D_u(k)}) and 
the second part is an application of Lemma \ref{approx_delta} below. 

Let $D_r(2)$ 
and $\widetilde{D}_r(2)$ be as defined above with the exception that  
the variables $x_1$ and $x_2$ are now integrated over $I_1$ and $I_v$ 
respectively. By construction $D_0(2)$ equals the left hand side in 
(\ref{D_{1,v}}). If we can show that 
\begin{equation} \label{D_r(2)}
\widetilde{D}_r(2) = \mathcal{O}(\sqrt{\epsilon})
\end{equation}
then by the same argument as above
\begin{equation*}
D_0(2) = D_{m-1}(2) + \mathcal{O}(\sqrt{\epsilon}).
\end{equation*}  
To see this we shall only need the trivial fact that
\begin{equation*}
|\fdet B| \leq \prod_{i=1}^n \sum_{j=1}^n |b_{ji}|
\end{equation*} 
where as before $B$ is a real $n \times n$ matrix.
Define $B$ as the $(m+2-r) \times (m+2-r)$ matrix appearing in 
$\widetilde{D}_r(2)$. We now estimate the column sums
\begin{equation*}
B_i := \sum_{j=1}^n |b_{ji}| .
\end{equation*} 

Column 1:
\begin{equation*}
B_1 = |A(p_{r+1},p_1)| + \ldots + |A(p_{m},p_1)| + 
|A(x_1,p_1)| + |A(x_2,p_1)| \leq Cm
\end{equation*}

Column 2:
\begin{multline*}
B_2 = \sqrt{\epsilon} \, \bigg(
|\widetilde{A}(p_{r+1},p_{r+2})| + |A(p_{r+2},p_{r+2})| + \ldots
+ |A(p_{m},p_{r+2})| \\ 
\qquad \qquad \qquad \qquad \qquad \qquad \qquad \qquad
+ |A(x_1,p_{r+2})| + |A(x_2,p_{r+2})| \bigg) \\
\leq \sqrt{\epsilon}  \,
\lp( Cm + |A(x_1,p_{r+2})| + |A(x_2,p_{r+2})| \rp)
\end{multline*}

Middle columns (if any) ($r+3 \leq i \leq m$):
\begin{equation*}
B_i = \sqrt{\epsilon} \, \lp( |A(p_{r+1},p_i)| + \ldots + |A(p_{m},p_i)| + 
|A(x_1,p_i)| + |A(x_2,p_i)| \rp)  \leq Cm
\end{equation*}

Last two columns:
\begin{multline*}
B_{m-r+1} = |A(p_{r+1},x_1)| + \ldots + |A(p_{m},x_1)| \\ + 
|A(x_1,x_1)| + |A(x_2,x_1)| \leq Cm e^{-x_1} 
\end{multline*}
\begin{multline*}
B_{m-r+2} = |A(p_{r+1},x_2)| + \ldots + |A(p_{m},x_2)| \\ + 
|A(x_1,x_2)| + |A(x_2,x_2)| \leq
C e^{-x_2} + \phi(x_1,x_2) + \sum_{k=r+1}^{v-1} \phi(p_k,x_2) 
\end{multline*}
Consider the estimates above for $B_2$ and $B_{m-r+2}$. 
The function $A(x_2,p_{r+2})$ will contain a $\phi$-function 
if and only if $v < r+2$, but in this case the sum
\begin{equation*}
\sum_{k=r+1}^{v-1} \phi(p_k,x_2)
\end{equation*}
is empty. This means that we do not get terms like 
\begin{equation*}
\phi(x_2,p_{r+2}) \phi(p_k,x_2)
\end{equation*}
in the product $B_2 B_{m-r+2}$. Given this
observation it is easy to see that
\begin{equation*}
\int_{I_1 \times I_v} B_2 B_{m-r+1} B_{m-r+2} \, \udx =
\mathcal{O}(\sqrt{\epsilon})
\end{equation*}
and this proves (\ref{D_r(2)}). 

The second part of the proof consists of showing that
\begin{multline} \label{sec_part}
\int_{I_1 \times I_v} \fdet \lp[
\begin{array}{ccc}
A(p_m,p_1) & A(p_m,x_1) & A(p_m,x_2) \\
A(x_1,p_1) & A(x_1,x_1) & A(x_1,x_2) \\
A(x_2,p_1) & A(x_2,x_1) & A(x_2,x_2) 
\end{array} \rp] \, \udx \\ =
\int_{I_1^2} R(p_1,x_1,x_2) \, \udx + 
\int_{I_1} R(p_1,x) \, \udx 
+ \mathcal{O}(\sqrt{\epsilon} \log{\epsilon}).
\end{multline}
The left hand side is equal to
\begin{multline*}
\int_{I_1 \times I_v} \fdet \lp[
\begin{array}{ccc}
A(p_m,p_1) & A(p_m,x_1) & A(p_m,x_2) \\
A(x_1,p_1) & A(x_1,x_1) & \widetilde{A}(x_1,x_2) \\
A(x_2,p_1) & A(x_2,x_1) & A(x_2,x_2) 
\end{array} \rp] \, \udx \\ +
\int_{I_1 \times I_v} \phi(x_1,x_2) \, \fdet \lp[
\begin{array}{cc}
A(p_m,p_1) & A(p_m,x_1) \\
A(x_2,p_1) & A(x_2,x_1) 
\end{array} \rp] \, \udx .
\end{multline*}
In view of Lemma \ref{K_approx} and (\ref{appr_delta_2}) 
in Lemma \ref{approx_delta} below we obtain (\ref{sec_part}).
\begin{lemma} \label{approx_delta}
Suppose that $f: \R \rightarrow \R$ has a continuous derivative and that
$g: \R^2 \rightarrow \R$ has continuous first partial derivatives. Assume that
\begin{equation*}
|f(x)|, \, |f'(x)| \leq C e^{-x}
\end{equation*}
\begin{equation*}
|g(x,y)|, \, |g'(x,y)| \leq C e^{-x-y}.
\end{equation*}  
Then, for $1 \leq i,j \leq m$, it holds that
\begin{equation}\label{appr_delta_1}
\int_{I_i} \frac{1}{\sqrt{4 \pi \epsilon}} 
e^{-\frac{(x-p_j)^2}{4 \epsilon}} f(x) \, \udx
= f(p_j) \, \int_{\frac{p_i-p_j}{\sqrt{\epsilon}}}^{\infty} 
\frac{1}{\sqrt{4 \pi}} 
e^{-\frac{x^2}{4}} \, \udx + \mathcal{O}(\sqrt{\epsilon})
\end{equation}
\begin{equation} \label{appr_delta_2}
\int_{I_i} \int_{I_j} 
\frac{1}{\sqrt{4 \pi \epsilon}} 
e^{-\frac{(x-y)^2}{4 \epsilon}} g(x,y) \, \udx \udy =
\int_{I_i} g(x,x) \, \udx + \mathcal{O}(\sqrt{\epsilon}\log{\epsilon}).
\end{equation}  
\end{lemma}
\begin{prf}
\begin{multline*}
\int_{p_i}^{\infty} \frac{1}{\sqrt{4 \pi \epsilon}} 
e^{-\frac{(x-p_j)^2}{4 \epsilon}} f(x) \, \udx = 
\lp[ z = \frac{x-p_j}{\sqrt{\epsilon}} \rp ] \\
= \int_{\frac{p_i-p_j}{\sqrt{\epsilon}}}^{\infty} \frac{1}{\sqrt{4 \pi}} 
e^{-\frac{z^2}{4}} f(p_j + \sqrt{\epsilon} z) \, \udz
\end{multline*}
By Taylors theorem
\begin{equation*}
f(p_j + \sqrt{\epsilon} z) 
= f(p_j) + \sqrt{\epsilon} z f'(p_j + \theta_{\epsilon}(z))
\end{equation*}
where $\theta_{\epsilon}(z)$ is a number between $0$ and $\sqrt{\epsilon}z$.
Since by assumption
\begin{equation*}
|f'(p_j + \theta_{\epsilon}(z))| \leq C e^{-p_j + \sqrt{\epsilon} |z|}
\end{equation*}
we obtain (\ref{appr_delta_1}).

\begin{multline*}
\int_{p_i}^{\infty} \int_{p_j}^{\infty} 
\frac{1}{\sqrt{4 \pi \epsilon}} 
e^{-\frac{(x-y)^2}{4 \epsilon}} g(x,y) \, \udx \udy
= \lp[ z = \frac{y-x}{\sqrt{\epsilon}} \rp ] \\
= \int_{p_i}^{\infty} \int_{\frac{p_j-x}{\sqrt{\epsilon}}}^{\infty} 
\frac{1}{\sqrt{4 \pi}} 
e^{-\frac{z^2}{4}} g(x,x + \sqrt{\epsilon} z) \, \udx \udz
\end{multline*}
By Taylors theorem
\begin{equation*}
g(x,x + \sqrt{\epsilon} z) 
= g(x,x) + \sqrt{\epsilon} z g'(x, x + \theta_{\epsilon}(x,z))
\end{equation*}
where $\theta_{\epsilon}(x,z)$ lies between $0$ and $\sqrt{\epsilon}z$.
The error can be discarded since
\begin{multline*}
\int_{p_i}^{\infty} \int_{\frac{p_j-x}{\sqrt{\epsilon}}}^{\infty} 
\frac{1}{\sqrt{4 \pi}} e^{-\frac{z^2}{4}} 
|z g'(x,x + \theta_{\epsilon}(x,z))| \, \udx \udz \\
\leq C \int_{p_i}^{\infty} \, \udx \int_{- \infty}^{\infty} 
\frac{1}{\sqrt{4 \pi}} |z| 
e^{-\frac{z^2}{4} - 2x + \sqrt{\epsilon} |z|} \, \udz
\leq C .
\end{multline*}
We now split the main term into two terms.
\begin{multline*}
\int_{p_i}^{\infty} \, \udx \int_{\frac{p_j-x}{\sqrt{\epsilon}}}^{\infty} \, \udz
\frac{1}{\sqrt{4 \pi}} e^{-\frac{z^2}{4}} g(x,x) \\
= \int_{p_i}^{p_i - \sqrt{\epsilon \log{\epsilon}}} \, \udx 
\int_{\frac{p_j-x}{\sqrt{\epsilon}}}^{\infty} \, \udz
\frac{1}{\sqrt{4 \pi}} e^{-\frac{z^2}{4}} g(x,x) \\ + 
\int_{p_i - \sqrt{\epsilon \log{\epsilon}}}^{\infty} \, \udx 
\int_{\frac{p_j-x}{\sqrt{\epsilon}}}^{\infty} \, \udz
\frac{1}{\sqrt{4 \pi}} e^{-\frac{z^2}{4}} g(x,x) =: \int_1 + \int_2
\end{multline*}
We can estimate the first integral by
\begin{equation*}
\lp| \int_1 \rp| \leq C \int_{p_i}^{p_i - \sqrt{\epsilon \log{\epsilon}}} \, \udx
\int_{-\infty}^{\infty} \, \udz \, e^{-\frac{z^2}{4} - 2x} \leq 
-C \sqrt{\epsilon} \log{\epsilon}.
\end{equation*}
If $x \geq p_j - \sqrt{\epsilon} \log{\epsilon}$ then 
$\frac{p_j-x}{\sqrt{\epsilon}} \leq C \log{\epsilon}$ and hence
\begin{multline*}
\int_{\frac{p_j - x}{\sqrt{\epsilon}}}^{\infty} 
\frac{1}{\sqrt{4 \pi}} e^{-\frac{z^2}{4}} \, \udz \\
= \int_{-\infty}^{\infty} \frac{1}{\sqrt{4 \pi}} 
e^{-\frac{z^2}{4}} \, \udz +
\int_{-\infty}^{\frac{p_j - x}{\sqrt{\epsilon}}} 
\frac{1}{\sqrt{4 \pi}} e^{-\frac{z^2}{4}}
= 1 + \mathcal{O}\lp( e^{- \frac{(\log{\epsilon})^2}{4}} \rp).
\end{multline*}
We finally get
\begin{multline*}
\int_2 = \int_{p_i - \sqrt{\epsilon} \log{\epsilon}}^{\infty}
\lp(  1 + \mathcal{O}\lp( e^{- \frac{(\log{\epsilon})^2}{4}} \rp) \rp)
g(x,x) \, \udx \\ = \int_{p_i}^{\infty} g(x,x) \, \udx +
\mathcal{O}(\sqrt{\epsilon} \log{\epsilon})
\end{multline*}
This concludes the proof of the lemma.
\end{prf}

\section{Theorem \ref{PNG}}

\subsection{Multi-layer discrete PNG}

Before we give the proof of Theorem \ref{PNG} we must present some preliminary 
results. 

How does one get a hand on the process $h$ described in the introduction? 
In \cite{Jo2} it is shown
that $h$ can be embedded as the top curve in a multi-layer process
given by a family of non-intersecting paths $\{ h_i, 0 \leq i < N \}$,
$h=h_0$. It turns out, see \cite{Jo2}, that
this multi-layer process is an example of a discrete determinantal
process. 
\begin{theorem}[Johansson]
Let $u,v \in \Z$ be such that $|u|,|v| < N$ and let $q = \alpha^2$. Set
\begin{equation*}
G(z,w) = (1 - \alpha)^{2(v-u)}
\frac{(1-\alpha/z)^{N+u} (1- \alpha w)^{N-v}}
{(1 - \alpha z)^{N-u} (1- \alpha / w)^{N+v}}
\end{equation*}
and 
\begin{equation*}
\widetilde{K}_N (2u,x;2v,y) =
\frac{1}{(2 \pi i)^2}
\int_{\gamma_{r_2}} \frac{\udz}{z} 
\int_{\gamma_{r_1}} \frac{\udw}{w} 
\frac{z}{z-w} G(z,w)
\end{equation*}
where $\gamma_r$ is the circle with radius $r$ centered around the origin,
$\alpha < r_1 < r_2 < 1/ \alpha$ and $x,y \in \Z$. Furthermore, define
\begin{equation*}
\phi_{2u,2v}(x,y) =
\frac{1}{2 \pi}
\int_{-\pi}^{\pi} e^{i(y-x) \theta} G \lp( e^{i \theta}, e^{i \theta}\rp) 
d \theta 
\end{equation*}
for $u<v$ and $\phi_{2u,2v}(x,y) = 0$ for $u \geq v$. Set 
\begin{equation*}
K_N (2u,x;2v,y) = \widetilde{K}_N(2u,x;2v,y) - \phi_{2u,2v}(x,y).
\end{equation*} 
Then,
\begin{multline*}
\Prob \big[ (2u,x_j^{2u}) \in \lp\{ 
(2t,h_i(2t,2N-1)); 0 \leq i < N , \, |t|<N \rp\},
\\ \qquad \qquad \qquad \qquad \qquad \qquad \qquad \qquad \qquad \qquad \qquad
|u| <N,1 \leq j \leq k_u \big] \\
= \fdet 
\lp( K_N(2u,x_i^{2u};2v,x_j^{2v}) \rp)_{|u|,|v|<N,1\leq i \leq k_u,1\leq j \leq k_v}
\end{multline*}
for any $x_j^{2u} \in \Z$ and any $k_u \in \{ 0,\ldots , N \}$.
\end{theorem} 
The asymptotic information about the kernel $K_N$ needed to prove Theorem
\ref{PNG} is contained in two lemmas. The first 
can be extracted from chapter four in \cite{Jo2} and the proof of the second
is provided at the end of this section. Please note that we make a slight
redefinition of the function $\phi$ from the last section. 
However, for the purposes of this text $\phi$ acts as one and the same.  
\begin{lemma}\label{Klemmat}
Let $\tau, \tau'$ be any real numbers such that
\begin{eqnarray*}
& & u = \frac{1 + \alpha}{1 - \alpha} d^{-1} N^{2/3} \tau \in \Z_+ \\
& & v = \frac{1 + \alpha}{1 - \alpha} d^{-1} N^{2/3} \tau' \in \Z_+.
\end{eqnarray*}
Let $x,y \in \Z_+$ and define $x',y'$ by
\begin{eqnarray*}
& & x = 2 \alpha (1 - \alpha)^{-1} N + (x' - \tau^2) d N^{1/3} \\
& & y = 2 \alpha (1 - \alpha)^{-1} N + (y' - \tau'^2) d N^{1/3} .
\end{eqnarray*}
For any $L \in \R$ there exist positive constants, $c$ and $C$, such that 
\begin{equation*}
| \widetilde{K}_N(2u,x;2v,y) | \leq C N^{-1/3} e^{- c(x'+y')} 
\end{equation*}
if $x',y' \geq L$.

If $|x'|,|y'| \leq \log{N}$, then the exists $c > 0$ such that
\begin{equation*}
d N^{1/3} \widetilde{K}_N(2u,x;2v,y) =
e^{\frac{\tau^3 - \tau'^3}{3} + y' \tau' - x' \tau} 
\widetilde{A}(\tau,x';\tau',y') + \Ordo(N^{-c}).
\end{equation*}
\end{lemma}
\begin{lemma}\label{philemmat}
Let $x,y \in \Z_+$ and define $x',y'$ by
\begin{eqnarray*}
& & x = 2 \alpha (1 - \alpha)^{-1} N + x' d N^{1/3} \\
& & y = 2 \alpha (1 - \alpha)^{-1} N + y' d N^{1/3} .
\end{eqnarray*}
Take $s > 0$, let $u \sim N^{2/3}$ and define $v$ by
\begin{equation*}
v = u + \frac{1 + \alpha}{1 - \alpha} d^{-1} s  N^{\gamma}
\end{equation*}
where $0 < \gamma < \frac{2}{3}$.
There exists a constant $C > 0$ such that
\begin{equation*}    
\phi_{2u,2v} (x,y) = \frac{1}{d N^{1/3}} \phi(x',y') + \phi_{E}(x',y')
\end{equation*}
where
\begin{equation*}
\phi(x',y') = \frac{1}{\sqrt{4 \pi s N^{\gamma - 2/3}}}
e^{-\frac{(x'-y')^2}{4 s N^{\gamma - 2/3}}}
\end{equation*}
and
\begin{equation*}
|\phi_{E}(x',y')| \leq \Bigg\{
\begin{array}{l}
C N^{-\frac{3 \gamma}{2}} \\
\frac{C}{N^{1/3} | x'-y'| N^{\gamma}}
\end{array}
\end{equation*}
for all $x,y$.
\end{lemma}

\subsection{Proof of Theorem \ref{PNG}}

This proof is really a discrete analog of the proof of Theorem \ref{theorem1}. 
Unfortunately things are more involved in this case where
$N^{\gamma - 2/3}$ plays the role of $\epsilon$. 

Please recall that $J_1 = \mu N + \psi d N^{1/3}$ where 
$\mu = 2 \alpha (1 - \alpha)^{-1}$ and $q = \alpha^2$.
Set $J_i = J_{i-1} + y_i d N^{\gamma / 2} \in \Z$, $i = 2, \ldots ,m$, and 
\begin{equation*}
\widetilde{I}_i = \lp\{ z \in \Z | z > J_i \rp\}.
\end{equation*}  
Here the $y_i$'s are arbitrary numbers such that $J_i \in \Z$. For later
convenience we also define $\psi_i$, $i = 1, \ldots ,m$, by 
$J_i = \mu N + \psi_i d N^{1/3}$.

We will prove that 
\begin{multline*}
\Prob \lp[ \# J_2 = \ldots = \# J_m =1, \# \widetilde{I}_2 = \ldots = 
\# \widetilde{I}_m = 0 | \# J_1 =1, \# \widetilde{I}_1 = 0 \rp] \\ =
\phi_{2K_1,2K_2}(J_1,J_2) \cdots \phi_{2K_{m-1},2K_m}(J_{m-1},J_m) \\ +
\Ordo \lp( ( N^{- \gamma /2} )^{m-1} N^{-c} \rp).
\end{multline*}
This implies Theorem \ref{PNG}: 
\begin{multline*}
\phi_{2K_1,2K_2}(J_1,J_2) \cdots \phi_{2K_{m-1},2K_m}(J_{m-1},J_m) \\ =
\frac{1}{\sqrt{4 \pi s_2}} e^{-\frac{y_2^2}{4s_2}} \cdots
\frac{1}{\sqrt{4 \pi s_m}} e^{-\frac{y_m^2}{4s_m}} 
\frac{1}{(d N^{\gamma / 2})^{m-1}} \lp( 1 + \Ordo \lp( N^{-c} \rp) \rp)
\end{multline*}
by Lemma \ref{philemmat}. The sum of this function over the sets $A_i$
is a Riemann sum that is well approximated by the integral in 
Theorem \ref{PNG}.

Define the finite integer intervals $I_i$, $1 \leq i \leq m$, by
\begin{equation*}
I_i = \lp\{ z \in \Z ; J_i < z < \lfloor \mu N \rfloor + N \rp\}.
\end{equation*}
The probability of finding a particle in $\widetilde{I}_i$ but outside of 
$I_i$ is very small:
\begin{align*}
\Prob & \lp[ \# ( \widetilde{I}_i \setminus I_i) \geq 1 \rp] \leq
\sum_{x \in \widetilde{I}_i \setminus I_i} \Prob [\# x = 1] =
\sum_{x \in \widetilde{I}_i \setminus I_i} K(x,x) \\ &  
= \sum_{k=0}^{\infty} 
K \Bigg(\lfloor \mu N \rfloor + \lp( \frac{1}{d}N^{2/3} + 
\frac{k}{dN^{1/3}} \rp)dN^{1/3}, \\ &
\qquad \qquad \qquad \qquad \qquad \qquad \qquad
\lfloor \mu N \rfloor + \lp( \frac{1}{d}N^{2/3} + 
\frac{k}{dN^{1/3}} \rp)dN^{1/3} \Bigg)
\\ & \qquad \qquad \qquad \qquad
\leq C \e ^{-\frac{1}{d} N^{2/3}} \sum_{k=0}^{\infty} e^{- \frac{k}{d N^{1/3}}} =
\Ordo \lp( e^{-c N^{2/3}} \rp).
\end{align*}
This means that we can work with $I_i$ instead of $\widetilde{I}_i$. We now
proceed much like we did in the proof of Theorem \ref{theorem1}.
If we set 
\begin{equation*}
A = \lp\{ \# J_1 = 1, \ldots, \# J_m = 1 \rp\}, 
\end{equation*}
then
\begin{multline*}
\Prob \lp[ A , \# I_1 = \ldots = \# I_m = 0 \rp] +
\Prob [ A , \# I_1 = 0, (\# I_2 = \ldots = \# I_m = 0)^c] \\ =
\Prob [ A , \# I_1 = 0]
\end{multline*}
where
\begin{multline*}
\Prob [ A , \# I_1 = 0, (\# I_2 = \ldots = \# I_m = 0)^c] \\ =
\Prob \lp[ A , \# I_1 = 0, \cup_{i=2}^m \{ \# I_i \neq 0 \} \rp] \\
\leq \sum_{i=2}^{m} \Prob [ A , \# I_1 = 0, \# I_i \neq 0 ] 
\leq \sum_{i=2}^{m} \Prob [ A , \# I_i \neq \# I_1 ],
\end{multline*}
and 
\begin{align*}
\Prob [ A , \# I_i \neq \# I_1 ]  & =
\E [\chi_{\{ \# J_1 = 1 \}} \cdots \chi_{\{ \# J_m = 1 \}} \cdot 
\chi_{\{ \# I_1 \neq \# I_i \}} ] \\
& = \E [ \# J_1 \cdots \# J_m \cdot \chi_{\{ \# I_1 \neq \# I_i \}} ] \\
& \leq \E [ \# J_1 \cdots \# J_m (\# I_1 - \# I_i )^2 ].
\end{align*}
The second equality holds since the probability of finding two particles
at the same place is zero. 

We need to prove three things: 
\begin{eqnarray*}
& & 1. \quad  \Prob [A, \# I_1 = 0] \\
& & \qquad \quad  
= \phi_{2K_1,2K_2}(J_1,J_2) \cdots \phi_{2K_{m-1},2K_m}(J_{m-1},J_m) \\
& & \qquad \qquad \qquad \qquad \qquad \qquad \qquad \qquad \qquad  \qquad
\times \Prob[ \# J_1 =1, \# I_1 = 0] \\
& & \qquad \qquad \qquad \qquad \qquad \qquad \qquad \qquad
+ \Ordo \lp( N^{-1/3-c} (N^{- \gamma / 2})^{m-1}\rp) \\
& & 2. \quad \E [ \# J_1 \cdots \# J_m \, (\# I_1 - \# I_i )^2 ] 
= \Ordo \lp( N^{-1/3-c} (N^{- \gamma / 2})^{m-1}\rp) \\
& & 3. \quad \Prob[ \# J_1 =1, \# I_1 = 0] \geq C N^{-1/3}
\end{eqnarray*}
Before giving the proofs we need some preliminaries.

When summing a function $f(x)$ over, say, $I_1$ we can write
\begin{equation*}
\sum_{x \in I_1}f(x) = \sum_{l=1}^{T_1} 
f\lp(\mu N + \lp( \psi_1 + \frac{l}{d N^{1/3}} \rp) d N^{1/3} \rp).
\end{equation*}
where $T_1 \sim N$. The next lemma will be frequently used later on.
\begin{lemma}\label{phi_sum}
There exists constants $C_1,C_2 > 0$ such that
\begin{equation*}
\sum_{k=1}^{\infty} \phi \lp( k / N^{1/3}, x \rp) N^{-1/3} \leq C_1
\end{equation*}
and 
\begin{equation*}
\sum_{k=1}^{N^2} \phi_E \lp( k / N^{1/3}, x \rp) \leq C_2 N^{-\gamma / 2}
\end{equation*}
for any $x \in \R$.
\end{lemma}
\begin{proof}
\begin{multline*}
\sum_{k=1}^{\infty} \phi \lp( k / N^{1/3}, x \rp) N^{-1/3} =
\sum_{k=1}^{\infty} \phi \lp( \frac{k - x N^{1/3}}{N^{1/3}}, 0 \rp) N^{-1/3} \\
\leq \sum_{k= - \infty}^{\infty} 
\phi \lp( \frac{k - x N^{1/3}}{N^{1/3}}, 0 \rp) N^{-1/3}
= \lp[ f := x N^{1/3} - \lfloor x N^{1/3} \rfloor \rp] \\ =
\sum_{k= - \infty}^{\infty} 
\phi \lp( \frac{k - f}{N^{1/3}}, 0 \rp) N^{-1/3} \\ 
\leq \sum_{k= - \infty}^{0} \phi \lp( \frac{k}{N^{1/3}}, 0 \rp) N^{-1/3} +
\phi \lp( \frac{1 - f}{N^{1/3}}, 0 \rp) N^{-1/3} +
\sum_{k=2}^{\infty} 
\phi \lp( \frac{k - 1}{N^{1/3}}, 0 \rp) N^{-1/3} \\ \leq
2 \sum_{k=1}^{\infty} 
\phi \lp( \frac{k}{N^{1/3}}, 0 \rp) N^{-1/3} + 2 \leq C
\end{multline*}
\begin{align*}
\sum_{k=1}^{N^2} & \phi_E \lp( k/N^{1/3},x \rp) \leq
C N^{-\gamma}\sum_{k=1}^{x N^{1/3} - N^{\gamma}} \frac{1}{x N^{1/3} - k} \\
& \qquad \qquad \qquad \qquad
+ C \sum_{x N^{1/3} - N^{\gamma}}^{x N^{1/3} + N^{\gamma}} N^{-3 \gamma/2} 
+ C N^{-\gamma}\sum_{x N^{1/3} + N^{\gamma}}^{N^2} \frac{1}{k - x N^{1/3}} \\
& \qquad \qquad
\leq C N^{-\gamma} \log{N} + C N^{-\gamma / 2} + C N^{-\gamma} \log{N} \leq
C N^{- \gamma / 2}
\end{align*}
\end{proof}
We now turn to the proof of 1. As in the proof of Theorem \ref{theorem1} we get
\begin{equation} \label{probeq}
\Prob [A, \# I_1 = 0] = \sum_{k=0}^{\infty} \frac{(-1)^k}{k!} 
\E \lp[ \# J_1 \cdots \# J_m \, \# I_1^{[k]} \rp] .
\end{equation}
For $0 \leq r \leq m-1$ set  
\begin{multline*}
D_r(k) \\ = \phi_{2K_1,2K_2}(J_1,J_2) \phi_{2K_2,2K_3}(J_2,J_3) \cdots 
\phi_{2K_r,2K_{r+1}}(J_r,J_{r+1}) 
\sum_{x_i \in I_1, 1 \leq i \leq k}  \\ \times \lp|
\begin{array}{ccccc}
K(J_{r+1},J_1) & K(J_{r+1},J_{r+2}) & \ldots & 
K(J_{r+1},J_{m}) & K(J_{r+1},x_j) \\
\vdots & \vdots & \, & \vdots & \vdots \\
K(J_{m},J_1) & K(J_{m},J_{r+2}) & \ldots & 
K(J_{m},J_{m}) & K(J_{m},x_j) \\
K(x_i,J_1) & K(x_i,J_{r+2}) & \ldots & 
K(x_i,J_{m}) & K(x_i,x_j) 
\end{array} \rp| .
\end{multline*}
The indicies $i,j$ run from $1$ to $k$ and if $r=0$ the (empty) product of
$\phi$-functions is to be interpreted as $1$. Let $\widetilde{D}_r(k)$ be
like $D_r(k)$ but having $\widetilde{K}(J_{r+1},J_{r+2})$ in position
(1,2) in the matrix. We want to show that
\begin{equation*}
|D_0(k) - D_r(k)| \leq  
N^{-1/3-c} \lp( N^{- \gamma/2} \rp)^{m-1} (Ck)^{\frac{k+m}{2}} 
\end{equation*}
which, by the induction argument in the proof of Theorem \ref{theorem1}, 
follows if we can prove that
\begin{equation} \label{pr_2_Dtilde}
|\widetilde{D}_r(k)| \leq 
N^{-1/3-c} \lp( N^{- \gamma/2} \rp)^{m-1} (Ck)^{\frac{k+m}{2}} .
\end{equation}
To show this we shall use Hadamard's inequality and therefore need to 
estimate sums of column elements squared 
(confer with the proof of Theorem 1). Lemmas \ref{Klemmat}, \ref{philemmat}
and \ref{phi_sum} will be frequently used below.

Column 1:
\begin{equation*}
\sum_{i=r+1}^m K^2(J_i,J_1) + \sum_{i=1}^k K^2(x_i,J_1) \leq
C N^{-2/3}(m+k)
\end{equation*}
\indent Column 2:
\begin{equation*}
\widetilde{K}^2(J_{r+1},J_{r+2}) + \sum_{i=r+2}^m K^2(J_i,J_{r+2}) \leq
C N^{-2/3} m
\end{equation*}
and
\begin{multline*}
\sum_{i=1}^k K^2(x_i,J_{r+2}) \\ \leq
C N^{-2/3} \sum_{i=1}^k 
\Big[ 1 + \phi \lp( l_i/dN^{1/3},\psi_1 - \psi_{r+2} \rp) \\ +
N^{1/3} \phi_E \lp(l_i/dN^{1/3},\psi_1 - \psi_{r+2} \rp) \Big]^2.
\end{multline*}

Columns $3, \ldots , m-r$ ($r+3 \leq j \leq m$), if they exist:
\begin{equation*}
\sum_{i=r+1}^m K^2(J_i,J_j) + \sum_{i=1}^k K^2(x_i, J_j) 
\leq C N^{-\gamma} (k+m)
\end{equation*}

Last $k$ columns ($1 \leq j \leq k$):
\begin{equation*}
\sum_{i=r+1}^m K^2(J_i,x_j) + \sum_{i=1}^k K^2(x_i, x_j) \leq
C (k+m) N^{-2/3} e^{- c l_j N^{-1/3}}
\end{equation*} 
Using Hadamard's inequality we get after some manipulations that
\begin{multline*}
|\widetilde{D}_r(k)| \\ \leq \sum_{l_1, \ldots ,l_k=1}^{T_1}
N^{-2/3} \lp(N^{-\gamma / 2}\rp)^{m-2}(Ck)^{\frac{k+m}{2}} \lp( N^{-1/3}\rp)^k 
\prod_{i=1}^k e^{- c l_i N^{-1/3}} \\ \times
\sum_{i=1}^k \Big[ 1 + \phi \lp(l_i/dN^{1/3},\psi_1 - \psi_{r+2} \rp)  \\ +
N^{1/3} \lp|\phi_E \lp(l_i/dN^{1/3},\psi_1 - \psi_{r+2} \rp) \rp| \Big] .
\end{multline*}
It follows from Lemma \ref{phi_sum} that
\begin{equation*}
\sum_{l_i=1}^{T_1} e^{-c l_i N^{-1/3}} 
\phi \lp(l_i/dN^{1/3},\psi_1 - \psi_{r+2} \rp) N^{-1/3} \leq C 
\end{equation*}
and also that
\begin{equation*}
\sum_{l_i=1}^{T_1} e^{-c l_i N^{-1/3}} 
\lp| \phi_E \lp(l_i/dN^{1/3},\psi_1 - \psi_{r+2} \rp) \rp|  \leq C N^{- \gamma / 2}
\end{equation*}
From this we get (\ref{pr_2_Dtilde}). 

To get 1 we also need to show that
\begin{multline}\label{det_sum}
\sum_{x_i \in I_1, 1 \leq i \leq k} 
\fdet \lp[
\begin{array}{cc}
K(J_m,J_1) & K(J_m,x_j) \\
K(x_i,J_1) & K(x_i,x_j)
\end{array} \rp]_{1 \leq i,j \leq k} \\
= \sum_{x_i \in I_1, 1 \leq i \leq k} 
\fdet \lp[
\begin{array}{cc}
K(J_1,J_1) & K(J_1,x_j) \\
K(x_i,J_1) & K(x_i,x_j)
\end{array} \rp]_{1 \leq i,j \leq k}
\\ + N^{-1/3 - c} \Ordo \lp( (Ck)^{\frac{k+m}{2}} \rp).
\end{multline}
Write 
\begin{equation*}
x_i = \mu N + \lp( \psi_1 + \frac{l_i}{d N^{1/3}} \rp)d N^{1/3}
\end{equation*}
and consider first the case $1 \leq l_i \leq N^{1/3} \log{N}$. From
Lemmas \ref{K_approx} and \ref{Klemmat} it is straight forward to deduce that
if $z=x_i$ or $z=J_1$ then
\begin{equation*}
K(J_m,z) = K(J_1,z) + \Ordo (N^{-1/3-c} ).
\end{equation*}
We now expand the determinant in the sum to the left in (\ref{det_sum}).
\begin{multline*}
\fdet \lp[
\begin{array}{cc}
K(J_m,J_1) & K(J_m,x_j) \\
K(x_i,J_1) & K(x_i,x_j)
\end{array} \rp]_{1 \leq i,j \leq k} \\
= \fdet \lp[ 
\begin{array}{cc}
K(J_1,J_1) & K(J_1,x_j) \\
K(x_i,J_1) & K(x_i,x_j)
\end{array} \rp]_{1 \leq i,j \leq k} \\
+ \Ordo (N^{-1/3-c} ) \sum_{p=1}^k
\fdet \lp[ K(x_i,J_1) \quad K(x_i,x_j) \rp]_{1 \leq i,j \leq k, j \neq p} \\
+ \Ordo ( N^{-1/3-c} )\, \fdet \lp[ K(x_i,x_j) \rp]_{1 \leq i,j \leq k}
\end{multline*}
We now use Hadamard's inequality to get
\begin{multline*}
\sum_{l_i=1}^{N^{1/3} \log{N}}
|\fdet \lp[ K(x_i,J_1) \quad K(x_i,x_j) \rp]_{1 \leq i,j \leq k, j \neq p}| \\
\leq   \sum_{l_i=1}^{N^{1/3} \log{N}}
(C k N^{-2/3} )^{k/2} 
e^{-N^{-1/3}(l_1 + \ldots + l_{p-1} + l_{p+1} + \ldots + l_k)} \\
\leq \lp( C k \rp)^{k/2} \log{N}
\end{multline*}
and
\begin{multline*}
\sum_{l_i=1}^{N^{1/3} \log{N}}
|\fdet \lp[ K(x_i,x_j\rp]_{1 \leq i,j \leq k}| \\
\leq \sum_{l_i=1}^{N^{1/3} \log{N}}
\lp(C k N^{-2/3} \rp)^{k/2} 
e^{-N^{-1/3}(l_1 + \ldots + l_k)} \leq \lp( C k \rp)^{k/2}.
\end{multline*}
This takes care of the summation over $1 \leq l_i \leq N^{1/3} \log{N}$, 
$1 \leq i \leq m$. By using Hadamard's inequality once more one readily
shows that the contribution coming from the remaining terms in the
sums in (\ref{det_sum}) is small enough to make (\ref{det_sum}) hold.

We now prove 2. Note that
\begin{equation*}
(\# I_1 - \#I_i)^2 = I_i^{[2]} + I_1^{[2]} + I_i + I_1 - 2 I_i  I_1.
\end{equation*}
By arguing as in the proof of 1 above we obtain
\begin{equation*}
\E \lp[ \# J_1 \cdots \# J_m \# I_u^{[k]} \rp]  =
\end{equation*}
\begin{multline*}
= \phi_{2K_1,2K_2}(J_1,J_2) \cdots \phi_{2K_{m-1},2K_m}(J_{m-1},J_m) \\
\qquad \qquad \qquad \qquad
\times \sum_{x_1, x_k \in I_u} 
\fdet \lp[ 
\begin{array}{cc}
K(J_m,J_1) &  K(J_m,x_s) \\
K(x_r,J_1) &  K(x_r,x_s) 
\end{array} \rp]_{1 \leq r,s \leq k}
\\ + \Ordo \lp( N^{-1/3-c - \frac{\gamma(m-1)}{2}} \rp)
\end{multline*}
where $u=1,i$ and $k=1,2$. One also gets
\begin{multline*}
\E [ \# J_1 \cdots \# J_m \# I_1 \# I_i]   \\
= \phi_{2K_1,2K_2}(J_1,J_2) \cdots \phi_{2K_{m-1},2K_m}(J_{m-1},J_m) \\
\times \sum_{x \in I_1, y \in I_i} 
\fdet \lp[ 
\begin{array}{ccc}
K(J_m,J_1) &  K(J_m,x) & K(J_m,y) \\
K(x,J_1) &  K(x,x) & K(x,y) \\
K(y,J_1) &  K(y,x) & K(y,y) 
\end{array} \rp]
\\ + \Ordo \lp( N^{-1/3-c - \frac{\gamma(m-1)}{2}} \rp).
\end{multline*} 
We omit the details. Using Lemma \ref{Klemmat} and Lemma
\ref{K_approx} one readily gets
\begin{multline*}
\sum_{x_1, x_k \in I_i} 
\fdet \lp[ 
\begin{array}{cc}
K(J_m,J_1) &  K(J_m,x_s) \\
K(x_r,J_1) &  K(x_r,x_s) 
\end{array} \rp]_{1 \leq r,s \leq k} \\
 = \E \lp[ \# J_1 \# I_1^{[k]} \rp] + \Ordo (N^{-1/3-c} )
\end{multline*}
for $k=1,2$ and
\begin{multline*}
\sum_{x \in I_1,y \in I_i} \fdet \lp[ 
\begin{array}{ccc}
K(J_m,J_1) &  K(J_m,x) & K(J_m,y) \\
K(x,J_1) &  K(x,x) & \widetilde{K}(x,y) \\
K(y,J_1) &  K(y,x) & K(y,y)
\end{array} \rp] \\
=  \E \lp[ \# J_1 \# I_1^{[2]} \rp] + \Ordo ( N^{-1/3-c} ).
\end{multline*}
We now see that 2 follows if 
\begin{multline} \label{phi_sum2}
\sum_{x \in I_1,y \in I_i} \phi_{2K_1,2K_i}(x,y)
\fdet \lp[ 
\begin{array}{cc}
K(J_m,J_1) &  K(J_m,x) \\
K(y,J_1) &  K(y,x)  
\end{array} \rp] \\
= \E [ \# J_1 \# I_1] + \Ordo (N^{-1/3-c} ).
\end{multline} 
We shall prove this by showing that both sides are well approximated by 
integrals. On the integral containing the function $\phi$ we can then apply 
Lemma \ref{approx_delta}.

By using Lemma \ref{phi_sum} we get rid of the error term associated 
with $\phi_E$:
\begin{multline*}
\sum_{l_1,l_2=1}^N 
\phi_E \lp(\psi_1 + \frac{l_1}{d N^{1/3}}, \psi_i + \frac{l_2}{d N^{1/3}} \rp)
e^{-\frac{l_1+l_2}{N^{1/3}}} N^{-2/3} \\
\leq  C \sum_{l_2=1}^N e^{-\frac{l_2}{N^{1/3}}} N^{-2/3-\gamma / 2}
\leq C N^{-1/3 - \gamma / 2}
\end{multline*}
The following calculation, again using Lemma \ref{phi_sum}, 
shows that the main contribution to the 
sums in (\ref{phi_sum2}) comes from summing over 
$1 \leq l_1,l_2 \leq N^{1/3} \log{N}$.
\begin{multline*}
\sum_{l_1=N^{1/3}\log{N}}^N \sum_{l_2=1}^N 
\phi \lp( l_2/dN^{1/3}, \psi_1 - \psi_i + l_2/dN^{1/3} \rp)
e^{-\frac{l_1+l_2}{N^{1/3}}} N^{-2/3} \\
\leq \sum_{l_1=N^{1/3}\log{N}}^N C e^{-l_1/N^{1/3}} N^{-1/3} \leq CN^{-1}
\end{multline*}
We shall use Euler's summation formula for two variables:
\begin{lemma}
Let $f(x,y)$ be a function of two variables such that its partial 
derivatives up to second order are continuous in the rectangle 
\begin{equation*}
\{ (x,y) |a \leq x \leq b, c \leq y \leq d \} 
\end{equation*}
where $a,b,c,d$ are integers. Then
\begin{align*}
\sum_{a \leq m \leq b} \sum_{c \leq n \leq d} f(m,n) = &
\int_a^b \int_c^d f(x,y) \, \udx \udy \\
& + \int_a^b \int_c^d f_x(x,y)(x - \lfloor x \rfloor) \, \udx \udy \\
& + \int_a^b \int_c^d f_y(x,y)(y - \lfloor y \rfloor) \, \udx \udy \\
& + \int_a^b \int_c^d f_{xy}(x,y)(x - \lfloor x \rfloor) 
(y - \lfloor y \rfloor)\, \udx \udy
\end{align*}
\end{lemma}
The case that we are interested in is when
\begin{equation*}
f(x,y) = 
\phi \lp( \psi_1 + \frac{x}{d N^{1/3}},\psi_i + \frac{y}{d N^{1/3}} \rp)
g \lp(x/d N^{1/3},y/d N^{1/3} \rp) N^{-1} 
\end{equation*}
where
\begin{equation*}
|g_x(x,y)|,|g_y(x,y)|, |g_x(x,y)| \leq  C e^{-c(x+y)} .
\end{equation*}
We need to show that the integrals involving the absolute values of
$f_x(x,y)$, $f_y(x,y)$ and $f_{xy}(x,y)$ are negligible. 
We only present the details for $|f_x(x,y)|$ here, the other terms 
are treated similarly.
\begin{multline*}
\int_1^{N^{1/3} \log{N}} \int_1^{N^{1/3} \log{N}} |f_x(x,y)| \, \udx \udy \\
\leq (d N^{1/3})^2 \int_{\psi_1}^{\infty} \int_{\psi_i}^{\infty}
\lp|f_x \lp((x-\psi_1) d N^{1/3},(y-\psi_i) d N^{1/3} \rp) \rp| \, \udx \udy \\
\leq C N^{-2/3} \int_{\psi_1}^{\infty} \int_{\psi_i}^{\infty}
\lp( |\phi_x(x,y)| + \phi(x,y) \rp) e^{-c(x+y)} \, \udx \udy
\end{multline*}
By Lemma (\ref{approx_delta}) 
\begin{equation*}
\int_{\psi_1}^{\infty} \int_{\psi_i}^{\infty}
\phi(x,y) e^{-c(x+y)} \, \udx \udy \leq C.
\end{equation*}
The remaining term demands some analysis.
\begin{align*}
\int_{\psi_1}^{\infty} & \int_{\psi_i}^{\infty}
|\phi_x(x,y)| e^{-c(x+y)} \, \udx \udy \\ & =
\int_{\psi_1}^{\infty} \int_{\psi_i}^{\infty} \frac{|x-y|}{2 N^{\gamma - 2/3}} 
\frac{1}{\sqrt{4 \pi N^{\gamma - 2/3}}}
e^{-\frac{(x-y)^2}{4 N^{\gamma - 2/3}}} e^{-c(x+y)} \, \udx \udy \\
& = \int_{\psi_1}^{\infty} \udx \lp( \int_{\psi_i}^x + \int_x^{\infty} \rp) 
\frac{|x-y|}{2 N^{\gamma - 2/3}} 
\frac{1}{\sqrt{4 \pi N^{\gamma - 2/3}}}
e^{-\frac{(x-y)^2}{4 N^{\gamma - 2/3}}} e^{-c(x+y)} \, \udy
\end{align*}
\begin{multline*}
\int_{\psi_1}^{\infty} \udx \int_{\psi_i}^x
\frac{x-y}{2 N^{\gamma - 2/3}} 
\frac{1}{\sqrt{4 \pi N^{\gamma - 2/3}}}
e^{-\frac{(x-y)^2}{4 N^{\gamma - 2/3}}} e^{-c(x+y)} \, \udy \\
= \int_{\psi_1}^{\infty} \udx \lp( 
\lp[ \phi(x,y) e^{-c(x+y)} \rp]_{\psi_i}^{x} +
c \int_{\psi_i}^x \phi(x,y) e^{-c(x+y)} \, \udy \rp) \\
\leq C N^{1/3 - \gamma / 2} + C \leq C N^{1/3 - \gamma / 2}
\end{multline*}
We can do the same calculation for the remaining integral. The 
$|f_x(x,y)|$ integral is hence $\Ordo \lp(N^{-1/3 - \gamma / 2}\rp)$ and
the same goes for the $|f_y(x,y)|$ and $|f_{xy}(x,y)|$ integrals.

Set
\begin{equation*} 
A^{\tau_1}(x,y) = A \lp(\tau_1,x + \tau_1^2;\tau_1, y + \tau_1^2 \rp).
\end{equation*}
Applying the above calculations to the left hand side of (\ref{phi_sum2})
and using Lemmas \ref{princ_lemma}, \ref{approx_delta}, \ref{K_approx} and 
\ref{Klemmat} we obtain
\begin{multline*}
\sum_{x \in I_1,y \in I_i} \phi_{2K_1,2K_i}(x,y)
\fdet \lp[ 
\begin{array}{cc}
K(J_m,J_1) &  K(J_m,x) \\
K(y,J_1) &  K(y,x)  
\end{array} \rp] \\
= \sum_{l_1,l_2 = 1}^{N^{1/3 \log{N}}}
\frac{1}{dN^{1/3}} \phi \lp( \psi_1 + l_1/ dN^{1/3},\psi_i + l_2/ dN^{1/3}\rp)
\lp( \frac{1}{dN^{1/3}} \rp)^2 \\
\times \lp|
\begin{array}{cc}
A^{\tau_1}(\psi_1, \psi_1) &  
e^{\tau_1 \frac{l_1}{dN^{1/3}}} 
A^{\tau_1}(\psi_1,\psi_1 + l_1/dN^{1/3}) \\
e^{-\tau_1 \frac{l_2}{dN^{1/3}}} 
A^{\tau_1}(\psi_1 + l_2/dN^{1/3},\psi_1) &  
A^{\tau_1}(\psi_i + l_2/dN^{1/3},\psi_1 + l_1/dN^{1/3})
\end{array} \rp| \\
\qquad \qquad \qquad \qquad \qquad \qquad \qquad \qquad
\qquad \qquad \qquad \qquad 
+ \Ordo(N^{-1/3-c})
\\= \frac{1}{d N^{1/3}} \int_{\psi_1}^{\infty}
\fdet \lp[ 
\begin{array}{cc}
A^{\tau_1}(\psi_1,\psi_1) &  A^{\tau_1}(\psi_1,x) \\
A^{\tau_1}(x,\psi_1) &  A^{\tau_1}(x,x)
\end{array} \rp] \udx
+ \Ordo(N^{-1/3-c}).
\end{multline*}
We get the same expression for the right hand side of (\ref{phi_sum2})
when applying Euler's summation formula. This concludes the proof of 2.

Let $F_2(t)$ be the Tracy-Widom distribution function corresponding to
the largest eigenvalue of the Gaussian Unitary Ensemble (GUE), \cite{TW}.
That 3 is true follows from the fact that $F_2'(t)>0$ $\forall t$, see 
\cite{TW}, together with the next lemma.
\begin{lemma} \label{TWconv}
Let $J_1$ and $I_1$ be as above. It holds that 
\begin{equation*}
\Prob [ \# J_1 = 0, \# I_1 = 0] = 
\frac{1}{d N^{1/3}} F_2'(\psi_1 + \tau_1^2) + \Ordo(N^{-2/3}).
\end{equation*}
\end{lemma}
\noindent \textbf{Proof:} This will, again, be an exercise in using 
Hadamard's inequality. We have the following representation for $F_2'$
(see the third equality in (\ref{twprimeq})):
\begin{equation}\label{TWdensity}
F_2'(t) = \sum_{k=0}^{\infty} \frac{(-1)^k}{k!} 
\int_{(t,\infty)^k} \fdet(A(x_i,x_j))_{0 \leq i,j \leq k} \, d^kx
\end{equation} 
where $x_0 = t$. In three steps we will now show that 
\begin{equation*}
d N^{1/3}\sum_{k=0}^{\infty} \frac{(-1)^k}{k!} 
\sum_{x_i \in I_1,1 \leq i \leq k} \fdet(K(x_i,x_j))_{0 \leq i,j \leq k}
\end{equation*}
where $x_0 = J_1$, is well approximated by the right hand side in 
(\ref{TWdensity}). By (\ref{probeq}) this will prove the lemma. In steps
one and two we will use Lemma \ref{Klemmat} to insert the kernel $A$ 
instead of $K$. In in the last step we show that we can change from summation
to integration. 

First we show that we can sum over 
$x_i = \mu N + \lp( \psi_1 + l_i/dN^{1/3} \rp)d N^{1/3}$ where
$1 \leq l_i \leq N^{1/3} \log{N}$, $1 \leq i \leq k$,
instead of over $I_1$. By Hadamard's inequality and Lemma \ref{Klemmat}
\begin{multline*}
\fdet \lp( K(x_i,x_j) \rp)_{0 \leq i,j \leq k} \leq 
\lp( \prod_{j=0}^k \sum_{i=0}^k K^2(x_i,x_j) \rp)^{1/2} \\
\leq \lp( C(k+1) N^{-2/3}\prod_{j=1}^k C(k+1) N^{-2/3}
e^{-c l_j / N^{1/3}} \rp)^{1/2} \\
\leq N^{-1/3} (C(k+1))^{\frac{k+1}{2}} \prod_{j=0}^k 
e^{-c l_j / N^{1/3}} N^{-1/3}.
\end{multline*}
We have that
\begin{multline*}
\sum_{\substack{l_i = 1 \\ 1 \leq i \leq k}}^{\infty}
\prod_{j=1}^k e^{-c l_j / N^{1/3}} N^{-1/3} -
\sum_{\substack{l_i = 1 \\ 1 \leq i \leq k}}^{N^{1/3} \log{N}}
\prod_{j=1}^k e^{-c l_j / N^{1/3}} N^{-1/3} \\
\leq k \sum_{l_1=1}^{N^{1/3} \log{N}} 
\sum_{\substack{l_i = 1 \\ 2 \leq i \leq k}}^{\infty}
\prod_{j=1}^k e^{-c l_j / N^{1/3}} N^{-1/3}
\leq k \, C^{k} N^{-1}.
\end{multline*}
Since 
\begin{equation*}
\sum_{k=1}^{\infty}\frac{1}{k!} N^{-1/3} (C(k+1))^{\frac{k+1}{2}}k N^{-1}
\leq C N^{-4/3}
\end{equation*}
we see that we can indeed restrict the summation.

In the second step we replace $K$ by $A$. As before we shall use the 
notation $A^{\tau}(x,y) = A(x+\tau^2,y+\tau^2)$. For 
$1 \leq l_i \leq N^{1/3} \log{N}$ it holds by Lemma \ref{Klemmat} that
\begin{multline*}
\fdet(K(x_i,x_j))_{0 \leq i,j \leq k} \\ =
\frac{1}{(dN^{1/3})^{k+1}} 
\fdet \lp(A^{\tau_1}(l_i/d N^{1/3},l_j/d N^{1/3}) + 
\Ordo(N^{-c})\rp)_{0 \leq i,j \leq k}
\end{multline*}
where we let $l_0 = \psi_1 d N^{1/3}$. If we expand the determinant in the right
hand side we get $(k+1)^2$ error terms of type
\begin{equation*}
\frac{N^{-c}}{(dN^{1/3})^{k+1}}  
\fdet \lp(A^{\tau_1}(l_i/d N^{1/3},l_j/d N^{1/3}) + 
\Ordo(N^{-c})\rp)_{\substack{0 \leq i,j \leq k \\ i \neq i_0, j \neq j_0}}.
\end{equation*}
An application of Hadamard's inequality together with Lemma \ref{Klemmat}
shows that the total error we get when changing from $K$ to $A^{\tau_1}$
is of order $N^{-1/3-c}$. We omit the details.

Finally we want to go from summation to integration. To do this we shall
use that
\begin{multline} \label{eulmacl1}
\sum_{l_i=1}^{N^{1/3} \log{N}}  
A^{\tau_1}(l_i /dN^{1/3},x)A^{\tau_1}(y,l_i/dN^{1/3}) \\
= dN^{1/3} \int_0^{\infty} A^{\tau_1}(z,x) A^{\tau_1}(y,z) \udz
+ \Ordo \lp(e^{-x-y} \rp) 
\end{multline}
and
\begin{multline}\label{eulmacl2}
\sum_{l_i=1}^{N^{1/3} \log{N}} A^{\tau_1}(l_i/dN^{1/3},l_i/dN^{1/3}) \\ =
dN^{1/3} \int_0^{\infty}A^{\tau_1}(z,z)\udz + \Ordo(1).
\end{multline}
This follows from Euler-Maclaurins summation formula and Lemma \ref{K_approx}.
We will show that
\begin{multline}\label{deteq}
\sum_{\substack{l_i=1 \\ 1 \leq i \leq k}}^{N^{1/3} \log{N}}
\frac{1}{(d N^{1/3})^{k+1}}
\fdet \lp( A^{\tau_1}(l_i/d N^{1/3},l_j/d N^{1/3}) \rp)_{0 \leq i,j \leq k} \\
= \frac{1}{d N^{1/3}} \int_{(0,\infty)^{k}}
\fdet \lp( A^{\tau_1}(y_i,y_j) \rp)_{0 \leq i,j \leq k} \, d^ky
\\ + \Ordo \lp( (Ck)^{\frac{k+5}{2}} N^{-2/3}\rp)
\end{multline}
where $l_0 = d N^{1/3} \psi_1$ and $y_0 = \psi_1$.
This will prove the lemma since 
\begin{equation*}
\sum_{k=1}^{\infty} \frac{1}{k!} (Ck)^{\frac{k+5}{2}} < \infty.
\end{equation*}
For $r = 0, \ldots ,k$ we set
\begin{equation*}
D_r = \frac{1}{(d N^{1/3})^{k-r+1}} 
\fdet \lp( A^{\tau_1}(z_i,z_j) \rp)_{0 \leq i,j \leq k}
\end{equation*}
where 
\begin{equation*}
z_i = \Bigg\{ 
\begin{array}{ll}
\psi_1 & \quad i=0 \\
y_i & \quad 1 \leq i \leq r \\
l_i/dN^{1/3} & \quad r+1 \leq i \leq k 
\end{array}. 
\end{equation*}
Please note that $D_0$ is what we sum over in (\ref{deteq}) and that $D_k$ is 
what we integrate over. $D_r$ should roughly be what we 
get after having changed summation over $l_1, \ldots ,l_r$ to integration over 
$y_1, \ldots ,y_r$. We can expand $D_r$ in such a way that we get
$k^2$ terms of type 
\begin{multline*}
\pm \frac{1}{(d N^{1/3})^{k-r+1}}
 A^{\tau_1}(z_{i_0},l_{r+1}/d N^{1/3}) A^{\tau_1}(l_{r+1}/d N^{1/3},z_{j_0}) \\
\times \fdet \lp( A^{\tau_1}(z_i,z_j)\rp)_{\substack{0 \leq i,j \leq k \\
i \neq r+1,i_0 \\ j \neq r+1,j_0}}
\end{multline*}
and one term 
\begin{equation*}
\frac{1}{(d N^{1/3})^{k-r+1}}
A^{\tau_1}(l_{r+1}/d N^{1/3},l_{r+1}/d N^{1/3})
\fdet \lp( A^{\tau_1}(z_i,z_j)\rp)_{\substack{0 \leq i,j \leq k \\
i,j \neq r+1}}.
\end{equation*}
We now apply (\ref{eulmacl1}) and (\ref{eulmacl2}) and therefore need to 
deal with the corresponding errors. 
\begin{multline*}
C (N^{1/3})^{k-r+1} e^{-z_{i_0} - z_{j_0}}
\fdet \lp( A^{\tau_1}(z_i,z_j) \rp)_{\substack{0 \leq i,j \leq k \\
i \neq i_0, r+1 \\ j \neq j_0, r+1}} \\
\leq C (N^{1/3})^{k-r+1} e^{-z_{i_0} - z_{j_0}}
\lp( \prod_{\substack{j=0 \\ j \neq j_0, r+1}}^k
C(k-1) e^{-c z_j} \rp)^{1/2} \\
\leq (N^{1/3})^{k-r+1} (C(k-1))^{\frac{k-1}{2}} 
\prod_{\substack{j=1 \\ j \neq r+1}}^k e^{-c z_j}
\end{multline*}
Since
\begin{equation*}
\int_{(0,\infty)^r} d^r x 
\sum_{\substack{l_i=1 \\ r+2 \leq i \leq k}}^{N^{1/3} \log{N}}
\prod_{\substack{j=1 \\ j \neq r+1}}^k e^{-c z_j}
\leq C^k (N^{1/3})^{k-(r+1)}
\end{equation*}
we find that the error from the $k^2$ terms of the first type is estimated by
\begin{equation*}
k^2 (C(k-1))^{\frac{k-1}{2}} N^{-2/3}.
\end{equation*}
The error coming from the remaining term can be treated in the same way. 
Changing
from summation over $l_i$ to integration over $y_i$, $1 \leq i \leq k$, 
hence results in an error estimated by
\begin{equation*}
k \, k^2 (C(k-1))^{\frac{k-1}{2}} N^{-2/3} = (Ck)^{\frac{k+5}{2}} N^{-2/3}
\end{equation*}
as needed.

\noindent \textbf{Proof of Lemma \ref{philemmat}:}
By definition
\begin{equation*}
\phi_{2u,2v} (x,y) = \frac{(1-\alpha)^{2(v-u)}}{2 \pi}
\int_{- \pi}^{\pi} e^{i(y-x) \theta + (u-v) 
\log{(1 + \alpha^2 - 2 \alpha \cos{\theta})}} d\theta.
\end{equation*}
Define
\begin{equation*}
g(\theta) = \log{(1 + \alpha^2 - 2 \alpha \cos{\theta})}
\end{equation*}
in $[-\pi,\pi]$.
This function is analytic in a neighbourhood of zero and a Maclaurin 
expansion gives
\begin{equation*}
g(\theta) = \log{(1-\alpha)^2} + \frac{\alpha}{(1-\alpha)^2 } \theta^2
+ c_2 \theta^4 + \mathcal{O}(\theta^6)
\end{equation*}
where $c_4 < 0$. It is easy to see that for any $\delta > 0$ there exists
$\epsilon > 0$ such that
\begin{equation*}
g(\theta) \geq \log{(1-\alpha)^2} + \epsilon
\end{equation*}
if $|\theta| \geq \delta$. Hence
\begin{equation*}
\lp| \int_{|\theta| > \delta} \frac{(1-\alpha)^{2(v-u)}}{2 \pi}
e^{i(y-x)\theta + (u-v)g(\theta)} d\theta \rp| 
\leq \frac{1}{2 \pi} 
\int_{\delta}^{\pi} e^{(u-v) \epsilon} d \theta
\sim e^{-\epsilon N^{\gamma}}.
\end{equation*}
We expect that the main contribution to $\phi_{2u,2v}$ will be
\begin{align*}
\frac{1}{2 \pi} \int_{-\delta}^{\delta}
& e^{i(y-x)\theta + (u-v) \frac{\alpha}{(1-\alpha)^2} \theta^2} 
\, d \theta   \\
& = \frac{1}{2 \pi} \int_{-\delta}^{\delta}
e^{i(y'-x') d N^{1/3} \theta - s d^2 N^{\gamma} \theta^2}
\, d \theta = \lp[ t = \sqrt{s}d N^{1/3}\theta \rp] \\
& = \frac{1}{2 \pi \sqrt{s} d N^{1/3}} 
\int_{-\delta \sqrt{s} d N^{1/3}}^{\delta \sqrt{s} d N^{1/3}}
e^{i\frac{y'-x'}{\sqrt{s}} t - N^{\gamma- \frac{2}{3}} t^2}
\, d t \\
& = \frac{1}{2 \pi \sqrt{s} d N^{1/3}} \int_{-\infty}^{\infty}
e^{i\frac{y'-x'}{\sqrt{s}} t - N^{\gamma- \frac{2}{3}} t^2}
\, d t + \mathcal{O}\lp( e^{-N^{\gamma}} \rp) \\
& = \frac{1}{d N^{1/3}} \frac{1}{\sqrt{4 \pi s N^{\gamma - 2/3}}}
e^{-\frac{(x'-y')^2}{4 s N^{\gamma - 2/3}}} + 
\mathcal{O}\lp( e^{-N^{\gamma}} \rp) .
\end{align*}
Below we will analyze the error. For simplicity we take $s=1$. 

Define $h(\theta)$ by
\begin{equation*}
g(\theta) = \log{(1-\alpha)^2} + \frac{\alpha}{(1-\alpha)^2}
\lp( \theta^2 + h(\theta) \rp).
\end{equation*}
This means that
\begin{equation*}
h(\theta) = \sum_{k=4}^{\infty} h_k \theta^k
\end{equation*}
where $h_4 < 0$. Note that $h$ is even since $g$ is and also that, 
for $\delta$ small enough, $h(\theta) < 0$ if $|\theta| \leq \delta$.
The error becomes
\begin{equation*}
\textrm{Err} = \lp| \int_{-\delta}^{\delta}
e^{i(y'-x')d N^{1/3} \theta} F(\theta) d \theta \rp| 
\end{equation*}
where
\begin{equation*}
F(\theta) = e^{- d^2 N^{\gamma} \theta^2} - 
e^{-d^2 N^{\gamma} \theta^2 - d^2 N^{\gamma}h(\theta)} .
\end{equation*}
Next we integrate by parts.
\begin{multline*}
\textrm{Err} \leq
\lp| \lp[ \frac{1}{i(y'-x')d N^{1/3}} 
e^{i(y'-x')d N^{1/3} \theta} F(\theta) \rp]_{-\delta}^{\delta} \rp| \\
\qquad \qquad + \frac{1}{|y'-x'|d N^{1/3}} \lp|
\int_{-\delta}^{\delta} 
e^{i(y'-x')d N^{1/3} \theta} F'(\theta) d \theta \rp| \\
\leq \frac{3}{|y'-x'|d N^{1/3}} e^{-d^2 N^{\gamma} \delta^2} + 
\frac{1}{|y'-x'|d N^{1/3}} \int_{-\delta}^{\delta}
\lp| F'(\theta) \rp| d \theta
\end{multline*}
The last integral will be easy to compute if we can find out where
$F'(\theta)$ changes sign. 
\begin{equation*}
F'(\theta) = 2 d^2 N^{\gamma} \theta \, e^{-d^2 N^{\gamma}(\theta^2 + h(\theta))}
\lp( 1 + \frac{h'(\theta)}{2 \theta} - e^{d^2 N^{\gamma}h(\theta)} \rp)
\end{equation*}
A point in $[-\delta,\delta] \setminus \{0\}$ where $F'$ changes sign 
will satisfy
\begin{equation*}
\frac{1}{d^2 N^{\gamma}} = 
\frac{h(\theta)}{\log{\lp[ 1 + \frac{h'(\theta)}{2 \theta} \rp]}} =
\frac{\theta^2}{2} + \Ordo(\theta^4).
\end{equation*}
This shows that if $N$ is large then $F'$ has two zeros $\pm \theta_0$ 
in $[-\delta,\delta] \setminus \{0\}$. Moreover, $\theta_0$ is of order 
$N^{-\gamma / 2}$.
Given this information we check which sign $F'$ has in different 
intervals and get
\begin{align*}
\int_{-\delta}^{\delta} |F'(\theta)| d \theta & =
2 \int_{0}^{\delta} |F'(\theta)| d \theta \\ & =
-\int_{0}^{\theta_0} F'(\theta) d \theta +
\int_{\theta_0}^{\delta} F'(\theta) d \theta \\ & =
F(0) - F(\theta_0) + F(\delta) - F(\theta_0) \\ & =
\Ordo (N^{-\gamma}).
\end{align*}
This almost finishes the proof of the second inequality in the lemma. 
We should not forget the exponentially small error terms that appeared 
above. They do not have the factor $|x'-y'|^{-1}$ in front of them.
However, a couple partial integrations can be used to take care of this
obstacle.    

The first inequality in the lemma follows from the following calculation.
\begin{multline*}
\int_{0}^{\delta} |F(\theta)| d \theta =
[ \theta = t N^{- \gamma / 2}] \\ =
N^{-\gamma / 2} \int_{0}^{N^{\gamma / 2} \delta}
e^{-d^2 t^2 - d^2 h(t N^{- \gamma / 2})} 
\lp( 1 - e^{d^2 N^{\gamma} h(t N^{- \gamma / 2})} \rp) d t \\ \leq
N^{-\gamma / 2} \int_{0}^{N^{\gamma / 2} \delta} e^{- c_1 t^2}
\lp( 1 - e^{-c_2 N^{-\gamma} t^4} \rp) d t \\ \leq
N^{-\gamma / 2} \int_{1}^{N^{\gamma / 2} \delta} t e^{- c_1 t^2}
\lp( 1 - e^{-c_2 N^{-\gamma} t^4} \rp) d t + \Ordo(N^{-3\gamma / 2})
\end{multline*}

We now use partial integration.
\begin{multline*}
\int_{1}^{N^{\gamma / 2} \delta} t e^{- c_1 t^2}
\lp( 1 - e^{-c_2 N^{-\gamma} t^4} \rp) d t \\ = 
\lp[ -\frac{1}{2 c_1} e^{- c_1 t^2}
\lp( 1 - e^{-c_2 N^{-\gamma} t^4} \rp) \rp]_1^{N^{\gamma / 2} \delta}  \\ +
\frac{2 c_2 N^{-\gamma}}{c_1} \int_{1}^{N^{\gamma / 2} \delta}
e^{- c_1 t^2} t^3 e^{-c_2 N^{-\gamma} t^4} d t \\
= \Ordo \lp( N^{-3 \gamma / 2} \rp)  
\end{multline*}
This concludes the calculations in this section as well as in this paper.

\section*{Acknowledgment}
I would like to express my deepest gratitude to my advisor Kurt Johansson
for his guidance and support.

\end{document}